\documentclass[a4paper,14]{elsarticle}

\usepackage[left=3cm,top=3cm,right=3cm,bottom=3cm]{geometry}
\usepackage[utf8]{inputenc}
\usepackage{times} %PostScript Times Roman instead of Computer Times Roman.
\usepackage{latexsym,amsfonts,amssymb,amsthm,amsmath,amscd,mathrsfs}
\usepackage{xargs}   % Use more than one optional parameter in a new commands
\usepackage{xspace}
\usepackage[shortlabels]{enumitem}
\usepackage{float}
\usepackage{multirow}
\usepackage[tight]{subfigure}
\usepackage{hyperref}
\usepackage{graphicx,tikz,times,pgffor}
\usepackage{pgfplots}
\usetikzlibrary{calc,through,backgrounds,patterns}
\usepackage{tabulary}
\usepackage{tabularx}
\usepackage{changepage}
\usepackage{soul}
\usepackage{float}

\pgfplotsset{compat=newest} 
\pgfplotsset{plot coordinates/math parser=false} 
\newlength\figureheight 
\newlength\figurewidth 

\usepackage{xpatch}
\usepackage{algpseudocode}
\usepackage{algorithm}

\newcommand{\setX}{\mathcal{X}}

\newcommand{\RN}{\mathcal{R}_0}
\newcommand{\Rn}{\mathcal{R}}
\newcommand{\Ko}{\mathcal{K}_0}
\newcommand{\K}{\mathcal{K}}
\newcommand{\Uc}{\mathcal{U}_c}

\newcommand{\X}{\mathbb{X}}

\newcommand {\bmat} {\left[\begin{array} }
\newcommand {\emat} {\end{array}\right]}

\newtheorem {theorem}{Theorem}
\newtheorem {defi}{Definition}
\newtheorem{rem}{Remark}
\newtheorem{lem}{Lemma}
\newtheorem {propt}{Property}

\newenvironment{pf}{\noindent \textbf{Proof:}}{}
\newenvironment{pf*}[1][Proof:]{\noindent \textbf{#1} }{}
\newenvironment{spf}{\noindent \textbf{Sketch of Proof:}}{}

\usepackage[shadow,colorinlistoftodos,prependcaption,spanish,textsize=footnotesize]{todonotes}
\numberwithin{equation}{section}
\numberwithin{theorem}{section}

\allowdisplaybreaks
\begin{document}
\begin{frontmatter}

\title{Dynamical Characterization of Antiviral Effects in COVID-19}

\author[First]{Pablo Abuin} 
\author[First]{Alejandro Anderson} 
\author[Second]{Antonio Ferramosca} 
\author[Third]{Esteban A. Hernandez-Vargas}
\author[First]{Alejandro H. Gonzalez} 

\address[First]{Institute of Technological Development for the Chemical Industry (INTEC), CONICET-UNL, Santa Fe, Argentina  (e-mail: alejgon@santafe-conicet.gov.ar)}
\address[Second]{Department of Management, Information and Production Engineering, University of Bergamo, Via Marconi 5, 24044, Dalmine (BG), Italy (e-mail: antonio.ferramosca@unibg.it)}
\address[Third]{Frankfurt Institute for Advanced Studies, 60438 Frankfurt am Main, Germany (e-mail: vargas@fias.uni-frankfurt.de)}

%\footnote{Corresponding author: alejgon@santafe-conicet.gov.ar, vargas@fias.uni-frankfurt.de}

%%%%%%%%%%%%%%%%%%%%%%%%%%%%%%%%%%%%%%%%%%%%%%%%%%%%%%%%%%%%%%%%%%%%%%%%%%%%%%%%
\begin{abstract}

Mathematical models describing SARS-CoV-2 dynamics and the corresponding immune responses in patients with COVID-19 can be critical to evaluate possible clinical outcomes of antiviral treatments. In this work, based on the concept of virus spreadability in the host, antiviral effectiveness thresholds are determined to establish whether or not a treatment will be able to clear the infection. In addition, the virus dynamic in the host  - including the time-to-peak and the final monotonically decreasing behavior - is chracterized as a function of the treatment initial time. Simulation results, based on nine real patient data, show the potential clinical benefits of a treatment classification according to patient critical parameters. This study is aimed at paving the way for the different antivirals being developed to tackle SARS-CoV-2.

\end{abstract}
\end{frontmatter}

%%%%%%%%%%%%%%%%%%%%%%%%%%%%%%%%%%%%%%%%%%%%%%%%%%%%%
\section{Introduction}\label{sec:intro}
%%%%%%%%%%%%%%%%%%%%%%%%%%%%%%%%%%%%%%%%%%%%%%%%%%%
With more than 82 million cases confirmed so far (December 2020) in 213 countries \cite{WHOCovid19SituationReport,Covid19Hopkins}, coronavirus disease COVID-19, caused by SARS-CoV-2 virus, continues spreading around the globe without neither effective treatment or vaccine strategy. In fact, the estimated case-fatality rate (CFR) for COVID-19 is about 5.7\%, while for the H1N1 pandemic the CFR was less than 1\% \cite{khan2020cardiovascular}.

Currently, a fully analysis about the potentiality of repurposed antiviral agents (\textit{i.e.}, Chloroquine and Hydroxychloroquine, Remdesivir, Favipiravir, Lopinavir/Ritonavir, Ribavirin) to ameliorate the viral spreading in the host are underway \cite{mitja2020use}.
The clinical observations suggest that prophylaxis with approved doses could prevent SARS-CoV-2 infection and reduce viral shedding \cite{Gao20,Yao20}, but these reports suffer from a number of limitations.
For instance, there is no certainty that the undefined benefits are not outweighed by the acute toxicity of the specific antiviral agents. To overcome this impasse, randomized clinical trials with adequate potency should be performed \cite{Kalil20}.
On the other hand, despite some reported dosing recommendation to treat COVID-19 \cite{Colson20,National20}, the efficacy of these proposed therapies cannot be guaranteed due to a paucity of data regarding the optimal dose of the antivirals \cite{sanders2020pharmacologic}.
In a race to find medical therapies to improve outcomes in patients with COVID-19 further studies are needed to elucidate the benefits of adapting antiviral agents. 

An important strategy to help to find the optimal dose of drugs is the pharmacological modelling based on in-vitro drug testing. This approach could suggest whenever the prophylaxis with an appropriate doses of antiviral agents could prevent SARS-CoV-2 infection or control the replication cycle of the virus \cite{Tett89}.
In this regard, recent studies have concentrated on the potential of a quantitative comprehension of COVID-19 dynamics \cite{vargas2020host,gonccalves2020timing,giordano2020sidarthe,kim2020modelling}.
Within-host mathematical models demonstrated useful insights about SARS-CoV-2 infection dynamics and its interactions with the immune system. More important, such results suggest the potential utility of assessing targets for drug development. 

The target cell-limited model, widely used to represent several diseases such as Influenza \cite{baccam2006kinetics,VargasJV14}, Ebola \cite{Nguyen15}, HIV \cite{Rong09,perelson2013modeling}, Hepatitis virus \cite{Graw16}, among others, has been linked to adjust the viral kinetics in infected patients with COVID-19 reported by \cite{wolfel2020virological}.
While a complete analysis of the main dynamic characteristic for the target cell model was developed to evaluate SARS-CoV-2 infection \cite{abuin2020char}, there is no adequate analysis of the effects of the existent pharmacological therapy for this model.

The main contribution of this work - that can be seen as an extension of \cite{abuin2020char} - is to provide a formal mathematical analysis of the COVID-19 dynamic under the effect of antiviral treatments, which may help to understand how to schedule the different therapies in function of the host parameters. A quantitative classification stating whether or not an antiviral will be effective - in terms of its capability to died out the virus in a reasonable period of time - is made, showing that their effectiveness could vary significantly between subjects.

After the introduction given in Section \ref{sec:intro} the article is organized as follows. 
Section \ref{sec:infectmod} proposes a general three-states "in host" model to represent COVID-19 infection dynamic and formalizes the concept of 
virus spreadability. 
Section \ref{sec:untreatsys} characterizes the equilibrium sets of the system, and its asymptotic stability, for the untreated case (\textit{i.e.}, when
no pharmacodynamic is considered), considering the system from the infection time $t=0$.
Then, in Section \ref{sec:pd}, the pharmacodynamic of antivirals is included, to analyze their effectiveness to avoid the virus to spread in the host.
Both, the inhibition effect on the replication rate of the virus and the inhibition effect on the infection rate of healthy cells are studied.
In Section \ref{sec:simres}, the results of Section \ref{sec:pd} are extensively demonstrated by simulating different treatments scenarios for nine patients from the literature.
Finally, Section \ref{sec:discuss} gives the discussion of the work, while several mathematical formalism - necessary to support the main results of article - are given in Appendix \ref{sec:app1}.
%%%%%%%%%%%%%%%%%%%%%%%%%%%%%%%%%%%%%%%%%%%%%%%%%%%%%%%%%%%%%%%%%%%%%%%%%%%%%%%%%%

%%%%%%%%%%%%%%%%%%%%%%%%%%%%%%%%%%%%%%%%%%%%%%%%%%%%%%%%%%%%%%%%%%%%%%%%%%%%%%%%%%
\section{Within-Host COVID-19 infection model}\label{sec:infectmod}
%%%%%%%%%%%%%%%%%%%%%%%%%%%%%%%%%%%%%%%%%%%%%%%%%%%

Mathematical models of within-host virus dynamic helps to improve the understanding of the interactions that
govern infections and permits the human intervention to moderate their effects \cite{hernandez2019modeling}. 
Basic models usually include the cells the pathogen infects, the pathogen particles, and their life 
cycle \cite{ciupe2017host} and, opposite to one can expect, they vary little in its structure from one infectious disease to another. 
Indeed, the well known target cell-limited models (\textit{i.e.}, models that incorporate the immune response just as a parameter) was successfully 
used to represent and control HIV \cite{perelson1993dynamics,legrand2003vivo,perelson2013modeling}, 
influenza \cite{larson1976influenza,baccam2006kinetics,smith2011influenza,hernandez2019passivity}, Ebola \cite{Nguyen15} and dengue \cite{nikin2015role,nikin2018modelling}, among others.
A main distinction between target cell models can be done according to the virus life cycle \cite{ciupe2017host}: chronic infection models (long-lived and persistent virus in comparison with cells life cycle) 
include production and death rates of healthy/susceptible cells, while acute infection models (short-lived virus) only consider the clearance of susceptible cells produced by the infection. This latter - which was firstly introduced in \cite{baccam2006kinetics} for the influenza - is used to describe the COVID-19, as it is detailed next.

%%%%%%%%%%%%%%%%%%%%%%%%%%%%%%%%%%%%%%%%%%%%%%%%%%%
\subsection{Target-cell-limited model} 
%%%%%%%%%%%%%%%%%%%%%%%%%%%%%%%%%%%%%%%%%%%%%%%%%%%

In this paper we consider the following mathematical model \cite{vargas2020host,abuin2020char}:
\begin{subequations}\label{eq:SysOrigAcut}
	\begin{align}
	&\dot{U}(t)   =  -\beta U(t) V(t),~~~~~~~ U(0) = U_0,\\
	&\dot{I}(t) =  \beta U(t) V(t) - \delta I(t),~~~~~~~I(0) = I_0=0,\\
	&\dot{V}(t)   = pI(t) - c V(t),~~~~~~~ V(0) = V_0,
	\end{align}
\end{subequations}
where $U$ $(cell)$, $I$ $(cell)$ and $V$ $(copies/mL)$ represent the uninfected cells, the infected cells, and the virus concentration.
The parameter $\beta$ $(mL day^{-1} /copies)$ is the infection rate, 
$\delta$ $(day^{-1})$ is the death rates of $I$,
$p$ $((copies)/(cell mL) day^{-1})$ is the viral replication, 
and $c$ $(day^{-1})$ is the viral clearance rate.  
The effects of immune responses are not explicitly described in this model, but they are implicitly included in the death rate of infected cells ($\delta$) and the clearance rate of virus ($c$) \cite{baccam2006kinetics}. %Figure \ref{figInf} shows a schematic plot of the interactions modeled in system \eqref{eq:SysOrigAcut}. 

System \eqref{eq:SysOrigAcut} is positive, which means that $U(t) \geq 0$, $I(t) \geq 0$ and $V(t) \geq 0$, for all $t\geq0$. We denote $x(t):=(U(t),I(t),V(t))$ the state vector, and
\begin{eqnarray}\label{eq:setX}
\mathbb X :=\{x\in \mathbb R^3_{\geq0}\},
\end{eqnarray}
the state constraints set. Another meaningful set (which is open) is the one consisting in all the states in $\mathbb{X}$ with strictly positive amount of virus and susceptible cells, \textit{i.e.},
\begin{eqnarray}\label{eq:setXcal}
\setX :=\{x\in \mathbb X : U>0,~V > 0\}.
\end{eqnarray}

The initial conditions of (\ref{eq:SysOrigAcut}), which represent a healthy steady state before the infection, are assumed to be $V(t) = 0$, $I(t) = 0$, and $U(t)= U_0>0$, for $t<0$. Then, at time $t=0$, a small quantity of virions enters the host body and, so, a discontinuity occurs in $V(t)$. Indeed, $V(t)$ jumps from $0$ to a small positive value $V_0$ at $t=0$ (formally, $V(t)$ has a discontinuity of the first kind at $t_0$, \textit{i.e.}, $\lim_{t\to0^-} V(t)=0$ while $\lim_{t\to0^+} V(t)=V_0>0$). 
This way, for the time after the infection, the virus may spread or being clear  depending on its infection effectiveness. To properly determine what such a spread means, the following (mathematical) definition is given
\begin{defi}[Spreadability of the virus in the host]\label{defi:spread}
	Consider system \eqref{eq:SysOrigAcut}, constrained by the positive set $\X$, with $V(t) = 0$, $I(t) = 0$, and $U(t)= U_0>0$, for $t<0$.
	Consider also that at time $t=0$, $V(t)$ jumps from $0$ to a small positive value $V_0$. Then, it is said that the virus spreads (in some degree) in the host for $t \geq 0$ if $\dot V(t) >0$ for some $t > 0$. If the virus does not spread in the host, it is said that it is cleared for $t \geq 0$.
\end{defi}
Definition \ref{defi:spread} states that the virus spreads in the host, for $t \geq 0$, if $V(t)$ increases at some time $t > 0$, and so, given that $\lim_{t\to\infty} V(t)=0$ (as it is stated later on in the article), it reaches at least one local maximum at some positive time. 
On the other hand, the virus is cleared for $t \geq 0$ if $V(t)$ is strictly decreasing for all $t>0$, which means that $V(t)$ has neither local minima nor local maxima at any $t >0$ (latter on it is stated that $V(t)$ has at most one minimum and one maximum). 
An infectious disease can be related to the virus peak and/or the period of permanence of the virus in the host \cite{Nguyen15}, and both effects are related to an increase of the virus load at some time after the infection. This is the reason why the spreadability is defined based on the virus positive derivative. A second reason supporting Definition \ref{defi:spread} is its use in the determination of antiviral effectiveness (Section \ref{sec:pd}). Antivirals able to clear the virus can be considered effective.
\begin{rem}
	Note that the virus spreadability may or may not imply a severe infection %an infection
	(a disease that eventually causes the host symptoms) depending on the value of the virus peak and the time the virus remains in the host. The non-spreadability of the virus, however, does imply that the infection is controlled. 
\end{rem}
In what follows, time $t=0$ will be considered as the infection time, \textit{i.e.}, the time when a small amount of virus enters the host. 
%

%%%%%%%%%%%%%%%%%%%%%%%%%%%%%%%%%%%%%%%%%%%%%%%%%%%
\subsection{Reproduction number and critical value for the target cells} 
%%%%%%%%%%%%%%%%%%%%%%%%%%%%%%%%%%%%%%%%%%%%%%%%%%%

To formally establish conditions under which the virus does or does not spread for $t>0$, some basic concepts need to be defined. The first one is the well-known basic reproduction number.
\begin{defi}\label{defi:repnumb}
	The within-host basic reproduction number $\Rn$ is defined as the number of infected cells (or virus particles) that are produced by one infected cell (or virus particle), over a course of its life-span. Its mathematical expression is given by:
	\begin{eqnarray}\label{eq:Rn}
		\Rn(t) := U(t)\frac{\beta p}{c \delta} .
	\end{eqnarray}
	Particularly, for $t=0$, this number describes the number of virus particles produced by a virus particle when a small amount of virus, $V_0$, is introduced into a healthy stationary population of uninfected target cells, $U_0$,
	\begin{eqnarray}\label{eq:RN}
		\RN := U_0\frac{\beta p}{c \delta} .
	\end{eqnarray}
\end{defi}
A second number, which is closely related to the first one, is the critical value for the susceptible cells such that the virus starts to decrease.   
\begin{defi}
	The critical number of susceptible cells $U$, $\Uc$, is defined as
	\begin{eqnarray}\label{eq:ucrit}
	\Uc := \frac{c \delta}{p \beta},
	\end{eqnarray}
	which for fixed system parameters $\beta$, $p$, $\delta$ and $c$ is a constant.
\end{defi}
Note that $U(t) < \Uc$ if and only if $\Rn(t) <1$, for every $t\geq0$. 
The basic reproduction number $\Rn(t)$ and the critical number $\Uc$ completely describe dynamic \eqref{eq:SysOrigAcut}, as it is shown in the next sections.
%%%%%%%%%%%%%%%%%%%%%%%%%%%%%%%%%%%%%%%%%%%%%%%%%%%%%%%%%%%%%%%%%%%%%%%%%%%%%%%%%%
%%%%%%%%%%%%%%%%%%%%%%%%%%%%%%%%%%%%%%%%%%%%%%%%%%%%%%%%%%%%%%%%%%%%%%%%%%%%%%%%%%
\section{Dynamical analysis for the untreated infection} \label{sec:untreatsys}
%%%%%%%%%%%%%%%%%%%%%%%%%%%%%%%%%%%%%%%%%
In this section, the infection evolution from the beginning (time $t=0$) is considered without  antiviral treatments. A detailed analysis of the dynamical behaviour will help us to properly understand, in next sections, the antivirals effects.
%%%%%%%%%%%%%%%%%%%%%%%%%%%%%%%%%%%%%%%%%%%%%%%%%%%
\subsection{Equilibrium set characterization and stability}\label{sec:equilset}
%%%%%%%%%%%%%%%%%%%%%%%%%%%%%%%%%%%%%%%%%%%%%%%%%%%
By equaling $\dot U$, $\dot I$ and $\dot V$ to zero, in \eqref{eq:SysOrigAcut}, it can be shown that the system only has free virus equilibrium of the form $x_s=(U_s,0,0)$, with $U_s$ being an arbitrary positive value, \textit{i.e.}, $U_s \in [0,\infty)$. Thus, there is only one equilibrium set, which is a healthy one, and it is defined by 
\begin{eqnarray}\label{eq:equilset}
\setX_s := \{(U,I,V)\in \mathbb{R}^3 : U \in [0,\infty),~I=0,~V=0\}.
\end{eqnarray}
According to \cite{abuin2020char}, set $\setX_s$ can be partitioned in two subsets, 
\begin{eqnarray}\label{ec:setXs1}
\setX_s^1 := \{(U,I,V)\in \mathbb{R}^3 : U \in [0,\Uc),~I=0,~V=0\}
\end{eqnarray}
and
\begin{eqnarray}\label{ec:setXs2}
\setX_s^2 := \{(U,I,V)\in \mathbb{R}^3 : U \in [\Uc,\infty),~I=0,~V=0\}
\end{eqnarray}
being $\setX_s^1$ asymptotically stable and $\setX_s^2$ unstable, as stated in the next theorem.
\begin{theorem}[Asymptotic stability]\label{theo:AS}
	Consider system \eqref{eq:SysOrigAcut} constrained by the positive set $\mathbb X$.
	Then, the set $\setX_s^1$ defined in \eqref{ec:setXs1} is the smallest asymptotically stable (AS) equilibrium set, with a domain of 
	attraction given by $\setX$.
\end{theorem}

\begin{pf}
	The proof is given in \cite{abuin2020char}, Theorem 3.3. $\square$
\end{pf}

\vspace{0.5cm}

Theorem \ref{theo:AS} establishes that neither single equilibrium points $x_s\in \setX_s^1$ nor subsets of $\setX_s^1$ are AS.
Furthermore, any trajectory starting at non-equilibrium states tends to a state in $\setX_s^1$, which depends on the initial conditions. This has critical implications on the controllability of the system, provided that $\setX_s^1$ is the minimal formal equilibrium target set to be used in a control strategy (based, for instance, on antiviral treatments).

Furthermore, as stated in \cite{abuin2020char}, if we define $U_\infty=\lim_{t\to\infty}U(t)$, $I_\infty=\lim_{t\to\infty}I(t)$ and $V_\infty=\lim_{t\to\infty}V(t)$, then, $I_\infty=0$, $V_\infty=0$, and
    \begin{eqnarray}\label{eq:Uinf}
    U_\infty=- \Uc W(-\RN e^{-(\RN+\Ko)}),
    \end{eqnarray}
where $W(\cdot)$ is the Lambert function, and $\RN = U_0 \frac{\beta p}{c \delta}$ and $\Ko = \frac{\beta}{c} (\frac{p}{\delta} I_0 - V_0)$ are constants depending on the initial conditions. Given that $W(-\RN e^{-(\RN+\Ko)}) \in (0,1)$ for positive values of $\RN$ and $\Ko$, this means that the only state that does not necessarily converge to zero, \textit{i.e.} $U$, converges to a value smaller than $\Uc$, which depends on $U_0,I_0,V_0$. See \cite{abuin2020char} for further details.

%%%%%%%%%%%%%%%%%%%%%%%%%%%%%%%%%%%%%%%%%%%%%%%%%%%
\subsection{System behavior} \label{sec:sysbeh}
%%%%%%%%%%%%%%%%%%%%%%%%%%%%%%%%%%%%%%%%%%%%%%%%%%%

The general behavior of system \eqref{eq:SysOrigAcut} can be analyzed according to the initial conditions and parameter values.
Property \ref{prop:infcond}, next, shows some general characteristics of the susceptible cells as time goes to infinity (\cite{abuin2020char}). 

\begin{propt} \label{prop:infcond}
	Consider system \eqref{eq:SysOrigAcut}, constrained by the positive set $\X$, at the beginning of the infection, \textit{i.e.},
	$U(0) = U_0 >0$, $I(0)= 0$ and $V(0)=V_0>0$ (\textit{i.e.}, $x(0)=(U(0),I(0),V(0)) \in \setX$). Consider also that $V_0$ is small enough. Then,
	\begin{enumerate}
		\item $U_\infty \to 0$ when $U_0\to\infty$ or $U_0\to 0$.
		\item $U_\infty \to \Uc$ when $U_0\to \Uc$.
		\item $0 < U_{\infty}( U_{0,1}, I_0, V_0) < U_{\infty}(U_{0,2},  I_0, V_0) < \Uc$, for initial conditions 
		$ U_{0,1} < U_{0,2} < \Uc$.
		\item $0 < U_{\infty}( U_{0,2}, I_0, V_0) < U_{\infty}( U_{0,1}, I_0, V_0) < \Uc$, for initial conditions
		$\Uc <  U_{0,1} < U_{0,2}$.
	\end{enumerate}
\end{propt}
\begin{pf}
	The proof - which is given in \cite{abuin2020char}, Property 2 - follows from the fact that $\Ko \approx 0$ for the selected initial conditions, and, so, according to \eqref{eq:Uinf}, $U_{\infty} \approx - \Uc W(-\RN e^{-\RN})$.
\end{pf}

\vspace{0.5cm}

\begin{figure}
	\centering
	\includegraphics[width=0.7\textwidth]{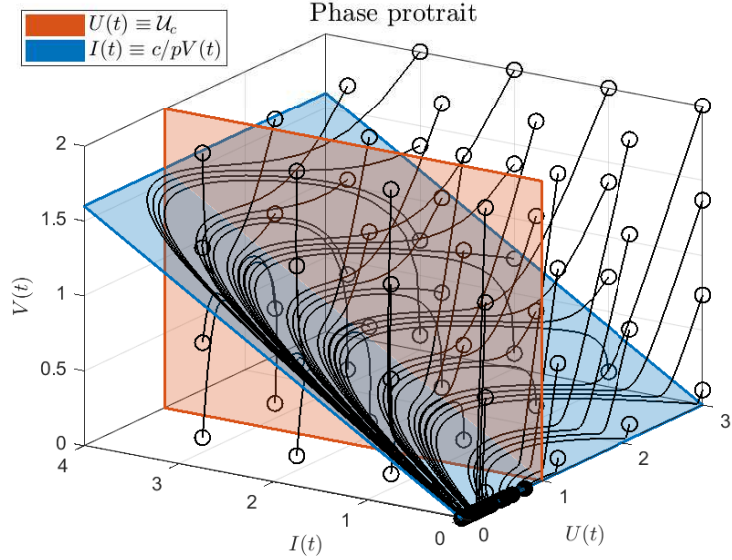}
	\caption{\small{Phase portrait of system \eqref{eq:SysOrigAcut} with parameters $\beta=0.5$, $\delta=0.2$, $p=2$ and $c=5$, for different initial conditions not necessarily representing realistic cases. Empty circles represent the initial state, while solid circles represent final states. The red hyperplane corresponds to $U(t) \equiv \Uc$ (\textit{i.e.}, the critical value of $U$, when $\Rn(t)=1$) while the blue hyperplane corresponds to the fast manifold in which $I(t)$ and $V(t)$ are proportional (.i.e, $I(t)=c/p V(t)$). Note that only the initial states with $U_0>\Uc=1$ corresponds to scenarios with $\RN>1$.}}
	\label{fig:PhasePortN}
\end{figure}

Figure \ref{fig:PhasePortN} shows a phase portrait of system \eqref{eq:SysOrigAcut}, with parameters $\beta=0.5$, $\delta=0.2$, $p=2$ and $c=5$, which provide $\Uc=1$. As it can be seen, every initial state, even those that are close to the axes $(U,0,0)$ with large values of $U$, tends to $\setX_s^1$ as $t \to \infty$.
The following theorem characterizes the whole dynamic of system \eqref{eq:SysOrigAcut} from the start of the infection, in terms of the variables peaks.

\begin{theorem}[Virus behavior from the infection time]\label{theo:key0}
	Consider system \eqref{eq:SysOrigAcut}, constrained by the positive set $\X$, at the beginning of the infection, \textit{i.e.},
	$U(0) = U_0 >0$, $I(0)= 0$ and $V(0)=V_0>0$ (\textit{i.e.}, $x(0)=(U(0),I(0),V(0)) \in \setX$). If the virus spreads (according to Definition~\ref{defi:spread}), 
	then $\Rn_0>1+\alpha(0)$, for some $\alpha(0)>0$ (or, the same, $U_0 > \Uc$) and there exist positive times $\check t_V$, $\hat t_I$, $t_c$ and $\hat t_V$, such that $\check t_V < \hat t_I < t_c < \hat t_V$, where $\check t_V$ and $\hat t_V$ are the times at which $V(t)$ reaches a local minimum and a local maximum, respectively, $\hat t_I$ is the time at which $I(t)$ reaches a local maximum, and $t_c$ is the time at which $U(t)$ reaches $\Uc$. Furthermore, $\dot V(t)<0$ for all $t> \hat t_V$.	
\end{theorem}

\begin{pf}
    The proof is given in \cite{abuin2020char}, Theorem 4.1.
\end{pf}

\begin{rem}\label{rem:alphaneg}
	Function $\alpha(t)$, $t\geq 0$, can be thought as a correction factor of the reproduction number $\Rn(t)$, necessary to understand and characterize the threshold over which the virus spreads in the host. As it is shown in \cite{abuin2020char}, Remark 5, this function cannot be explicitly defined but may be computed numerically. In general, for real patient data, it takes small values when the virus load is small (particularly, at time $t=0$), and becomes significant only when the virus approaches its peak, at $\hat t_V$. Given that in most of the interesting situations (\textit{i.e.}, at the beginning of the infection or when an antiviral is administered before the virus peak) its effect is negligible, $\alpha(t)$ will not be considered in the computation of critical values for the treatments. 
\end{rem} 
\begin{rem}\label{rem:sufcond}
	Note that from Theorem \ref{theo:key0}, it follows that a sufficient condition (not necessary) for the virus not to spread in the host  at time $t>0$ is given by $\RN<1$ (or, the same, $U(0)<\Uc$). This condition, adapted for the time an antiviral treatments is initiated, will serve in the next section to determine their effectiveness.
\end{rem}

Next, a Remark concerning a particularity of Theorem \ref{theo:key0} is introduced, which may help to approximately determine when the maximum of the virus load occurs.

\begin{figure}
	\centering
	\includegraphics[width=0.8\textwidth]{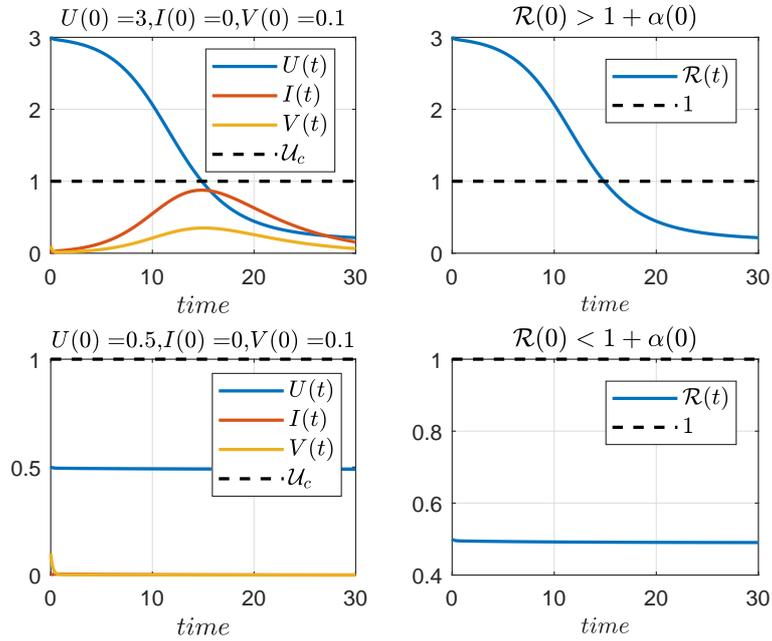}
	\caption{\small{Time evolution of $U(t)$, $I(t)$, $V(t)$, and $\Rn(t)$ corresponding to system \eqref{eq:SysOrigAcut}, with parameters $\beta=0.5$, $\delta=0.2$, $p=2$ and $c=5$ and initial conditions $U_0=3,I_0=0,V_0=0.1$ (first row) and $U_0=0.5,I_0=0,V_0=0.1$ (second row). The firs row corresponds to the case $\RN >1+\alpha(0)$, while the second one corresponds to the case $\RN<1+\alpha(0)$ (with $\alpha(0)$ being approximately zero). In the first case, the minimum of $V(t)$ occurs at $\check t_V=0.85$, the maximum of $I$ occurs at $\hat t_I=14.86$, $U(t)$ crosses $1$ at $t_c=14.87$ and the maximum of $V(t)$ occurs at $\hat t_V=15.08$, which confirms both, the results in Theorem \ref{theo:mainresult} and Remark \ref{cor:Vmax}. In the second case, both $I(t)$ and $V(t)$ are strictly decreasing and $U(t)$ remains almost constants, which confirms the results in Theorem \ref{theo:mainresult} and Remark \ref{rem:sufcond}.}}
	\label{fig:main}
\end{figure}

\begin{rem}\label{cor:Vmax}
	Consider that the assumptions of Theorem \ref{theo:key0} hold. If the virus clearance is significantly faster than the infected cell death ($c>>\delta$), as it is usually the case (see \cite{kim2020modelling}, Equation (4), \cite{ikeda2016quantifying}, \cite{nowak2000virus}, \cite{dahari2007modeling}) the three-states system \eqref{eq:SysOrigAcut} can be approximated by the following two-state equations:
	\begin{subequations}\label{eq:SysOrigAcutApp}
		\begin{align}
		&\dot{U}(t)   \approx  -\beta U(t) V(t),~~~~~~~ U(0) = U_0,\\
		&\dot{V}(t) \approx  \frac{p \beta}{c} U(t)V(t) - \delta V(t),~~~~~~~V(0) = V_0,
		\end{align}
	\end{subequations}
	where the infected cell state is given by $I(t) \approx \frac{c}{p}V(t)$. 
	Note that equation (\ref{eq:SysOrigAcutApp}.b) can be written as $\dot{V}(t)=(\Rn(t)-1) \delta V(t)$. Then, according to Theorem \ref{theo:key0},
	it is easy to see that $\hat t_V \to t_c$ from the right and $\hat t_I \to t_c$ from the left when $(c/\delta)  \to \infty$, \textit{i.e.}, the peaks of $V$ and $I$ tend to occur simultaneously, at time $t_c$. Figure \ref{fig:PhasePortN} shows a phase portrait of system \eqref{eq:SysOrigAcut} with a rather unrealistic parameter values ( $\beta=0.5$, $\delta=0.2$, $p=2$ and $c=5$) but useful to exemplify how, after a relatively short time (and provided that $c>\delta$), the system reaches the manifold defined by the condition $I(t) = \frac{c}{p}V(t)$. This behavior can also be seen in the time-behavior plot in Figure \ref{fig:main}.
\end{rem}

Figure \ref{fig:main} shows $U(t)$, $I(t)$, $V(t)$ and $\Rn(t)$ time evolution, for $\beta=0.5$, $\delta=0.2$, $p=2$ and $c=5$, and initial conditions such that $\RN >1+\alpha(0)$ and $\RN<1+\alpha(0)$ (with $\alpha(0) \approx 0$). As it can be seen, the results in Theorem \ref{theo:key0} and Remarks \ref{cor:Vmax} and \ref{rem:sufcond} are confirmed by the simulations.
All these results are also confirmed by the nine real patients identified and simulated in Section \ref{sec:pac} and \ref{sec:simres} (see, particularly, Table \ref{tab:const} and Figure \ref{fig:allpatOL}).

%%%%%%%%%%%%%%%%%%%%%%%%%%%%%%%%%%%%%%%%%%%%%%%%%%%%%%%%%%%%%%%%%%%%%%%%%%%%%%%%%%

%%%%%%%%%%%%%%%%%%%%%%%%%%%%%%%%%%%%%%%%%%%%%%%%%%%%%%%%%%%%%%%%%%%%%%%%%%%%%%%%%%
\section{Antiviral Treatments} \label{sec:pd}
%%%%%%%%%%%%%%%%%%%%%%%%%%%%%%%%%%%%%%%%%%%%%%%%%%%%%

The goal of this section is to formally consider the effect of antivirals into system \eqref{eq:SysOrigAcut} to obtain a controlled system, \textit{i.e.} a system with certain control actions - given by the antivirals - that allows one to (even partially) modified the whole system dynamic according to some control objectives.
Antivirals have the potential to inhibit the virus replication,  reducing the advance of the infection over the target cells of infected host (\textit{i.e.}: epithelial cells in the  respiratory tract for H1N1 \cite{baccam2006kinetics}).
Several antivirals are being tested for COVID-19 treatment, with different results concerning both, their inhibition effect on the virus replication and their toxicity. Among them, Remdesivir, Favipiravir, Umifenovir, Chloroquine, Oseltamivir, etc. \cite{jomah2020clinical, sanders2020pharmacologic,liu2020research,wang2020remdesivir,dobrovolny2020quantifying} can be mentioned.
The antiviral effect can be modeled as a reduction of the virus infectivity in the presence of inhibitors (by reducing the infection rate $\beta$) and/or as a reduction in the replication of infectious virions (by reducing the replication rate $p$) \cite{mitja2020use,kim2020modelling}.
In any case, the effectiveness of a treatment is limited and depends on the patient parameters (which in turn depends 
on his/her clinical state). Next, two different inhibition effects $\eta_{\beta}(t)$ and $\eta_p(t)$ are considered,
one affecting the infection rate $\beta$ and the other affecting the replication rate $p$. It is assumed that the virus spreads in the host from the infection time $t=0$, \textit{i.e.}, $V(t)$ reaches a maximum at some time $\hat t_V$ if no treatment is initiated. The controlled system (\textit{i.e.}, the system considering antiviral treatment) can be written as:
\begin{subequations}\label{eq:SysOrigAcutAntiv}
	\begin{align}
	&\dot{U}(t)   =  - \beta(1 - \eta_{\beta}(t)) U(t) V(t),  &U(0)=U_0, \\
	&\dot{I}(t) =  \beta(1 - \eta_{\beta}(t)) U(t) V(t) - \delta I(t),  &I(0)=I_0, \\
	&\dot{V}(t)   = p (1 - \eta_{p}(t)) I(t) - c V(t),  &V(0)=V_0,
	\end{align}
\end{subequations}
where $\eta_{\beta}(t) \in [0,1)$ and $\eta_p(t)\in [0,1)$ are assumed to jump from $0$ to the new value $\eta_{\beta}$ and $\eta_p$ at the treatment time $t_{tr}$ (\textit{i.e.}, the pharmacokinetic of the antivirals is assumed to have time constants significantly smaller than the ones of system \eqref{eq:SysOrigAcut}, and, so, it is neglected)
\begin{eqnarray}\label{eq:etas}
	\eta_{\beta}(t) = \left\lbrace 
	\begin{array}{cc}
		0   &    t < t_{tr}   \\
		\eta_{\beta}   &  t\geq t_{tr}  
	\end{array} \right., ~~~~
		\eta_{p}(t) = \left\lbrace 
	\begin{array}{cc}
	0   &    t < t_{tr}   \\
	\eta_{p}   &  t\geq t_{tr}  
	\end{array} ,\right. 
\end{eqnarray}
being $\eta_{\beta}\in [0,1)$ and $\eta_p\in [0,1)$ constant values representing the full inhibition treatment effect. The treatment time $t_{tr}$ is assumed to be between the minimum and maximum time of $V$, \textit{i.e.}, $\check t_V < t_{tr} < \hat t_V$ (although some simulations are performed for $t_{tr} > \hat t_V$).  
The full antiviral effect $\eta_p$ and $\eta_{\beta}$ are limited by the inhibitory potential of the drug (expressed in terms of EC$_{50}$, or drug concentration for inhibiting $50\%$ of antigen particles) and its cytotoxic effect (expressed in terms of IC$_{50}$, or drug concentration which causes death to $50\%$ of susceptible cells) \cite{wang2020remdesivir,vergnaud2005assessing}. 

As the antiviral treatment reduces the system parameter in some amount, it will quantitatively modify the virus behavior. Particularly, the virus peak time will be modified from $\hat t_V$ (untreated patient case) to $\hat t_{V,tr}$ (treated patient case). However, given that the treatment is initiated when the virus is increasing (\textit{i.e.}, between $\check t_V$ and $\hat t_V$), then the new peak will occur at the same time or after the treatment time, \textit{i.e.}, $\hat t_{V,tr} \geq t_{tr}$. This way, even when the virus peak will always be smaller with a treatment (smaller peaks are obtained for smaller values of $p$ and $\beta$), the time at which this peak takes place can be smaller or greater than the one without any treatment. This effect, usually disregarded in many studies concerning the effectiveness of antivirals, could be critical to define whether or not a given antiviral is able to impede a severe disease. Indeed, in some cases, antivirals significantly delay the virus peaks, largely increasing the time of permanence of the virus in the host. In order to qualitatively assess antiviral effectiveness according to the time of the virus peak, the following classification is made:
\begin{defi}[Antiviral treatment effectiveness]\label{def:anttreat}
Consider system \eqref{eq:SysOrigAcutAntiv}, constrained by the positive set $\X$, such that the virus spreads in the host from time $t=0$, which implies that: $U(0) > \Uc$, $I(0) = 0$ and $V(0) > 0$. Consider also that, at time $t_{tr}$, with $\check t_V < t_{tr} < \hat t_V $, an antiviral treatment is initiated such that $\eta_p(t)$ and/or $\eta_{\beta}(t)$ jump from $0$ to $\eta_p \in [0,~1)$ and/or $\eta_{\beta} \in [0,~1)$, respectively (as stated in \eqref{eq:etas}). Then, the treatment is said to be effective if the virus peaks at a time $\hat t_{V,tr}<\hat t_V$, being the latter the virus peak time for the untreated virus evolution (\textit{i.e.}, when $\eta_p=\eta_{\beta}=0$). Otherwise, if $\hat t_{V,tr} \geq \hat t_V$, it is said that the treatment is ineffective.
\end{defi}
Definition \ref{def:anttreat} is closely related to the capacity of the antiviral drug to clear the virus infection (or, the same, cutting off its spread) in such a way that it could: a) decline the viral grow at the treatment time, and, so, the virus eradication starts when the therapy is initiated, or b) hasten the virus peak, and, so, even though the virus eradication is not started at treatment time, it begins prior to the untreated case. Note also that Definition \ref{def:anttreat} accounts for three typical antiviral effect metrics: the area under the virus curve (AUC) and the duration of infection (DI) \cite{hadjichrysanthou2016understanding}, in a direct way, and the difference of viral loads at the time-to-peak $\Delta V$\cite{gonccalves2020timing}, in an indirect way.
%
%%%%%%%%%%%%%%%%%%%%%%%%%%%%%%%%%%%%%%%%%%%%%%%%%%%
\subsection{Antiviral effectiveness characterization}\label{sec:andcharact}
%%%%%%%%%%%%%%%%%%%%%%%%%%%%%%%%%%%%%%%%%%%%%%%%%%%

\begin{figure}[hbt!]
	\centering
	\includegraphics[width=0.7\textwidth]{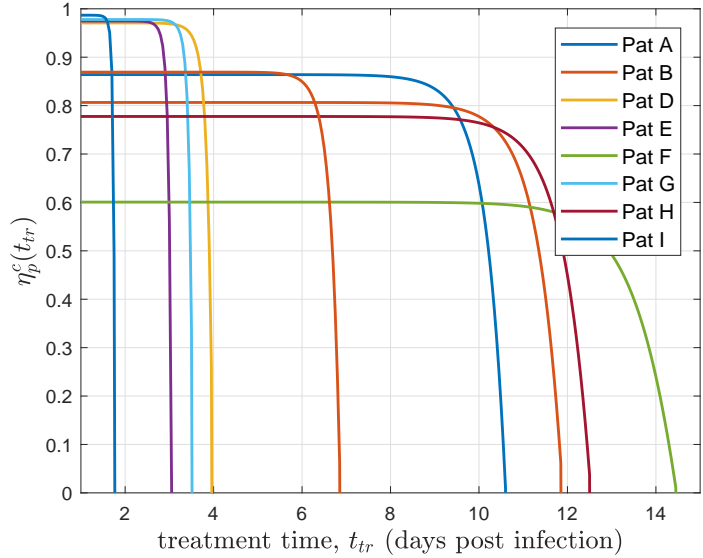}
	\caption{$\eta_p^c(t_{tr})$ vs $t_{tr}$ corresponding to the nine patients simulated in Section \ref{sec:pac}}
	\label{fig:eta_c_tr}
	\centering
\end{figure}

In this subsection it is shown that the effectiveness of antivirals depends on weather $\eta_p$ and/or $\eta_{\beta}$ are greater or smaller than a specified threshold, which is a function of the parameters and the time of the treatment initiation. In order to characterize such thresholds, the within-host basic reproduction number $\Rn$ at treatment time $t_{tr}$ is computed as follows: 
\begin{eqnarray}\label{eq:Rtreat}
\Rn(t_{tr})= \frac{U(t_{tr}) p(1-\eta_p) \beta(1-\eta_{\beta})}{c \delta},
\end{eqnarray} 
where $\check t_V < t_{tr} < \hat t_V$. The critical values of $\eta_p$ and $\eta_{\beta}$ are the ones that make $\Rn(t_{tr}) =1$, \textit{i.e.}:
\begin{eqnarray}
\eta_p^c(t_{tr})&:=& 1 -  \frac{c \delta}{U(t_{tr}) \beta p}\label{eq:eta_p}\\
\eta_{\beta}^c(t_{tr})&:=& 1 -  \frac{c \delta}{U(t_{tr}) \beta p} = \eta_p^c(t_{tr}) \label{eq:eta_beta}.
\end{eqnarray}
From equations \ref{eq:eta_p} and \ref{eq:eta_beta} it can be inferred that $\eta_p^c(t_{tr})$ and $\eta_\beta^c(t_{tr})$ are increasing functions of $U(t_{tr})$. Figure \ref{fig:eta_c_tr} shows the time behavior of $\eta_p^c(t_{tr})$ for the nine COVID-19 patients identified in Section \ref{sec:pac}. Note that $\eta_p^c(t_{tr})\approx 1-c\delta/(U_0\beta p)$ for $t_{tr}\rightarrow t_0$ and $\eta_p^c(t_{tr})\approx 1-c\delta/(\Uc\beta p)=0$ for $t_{tr}\rightarrow \hat t_V$ ($t_c \approx \hat t_V$). Similar results concerning the critical drug efficacy with respect to the availability of target cells at the treatment time were reached in \cite{dobrovolny2011neuraminidase}, for acute models, although the authors focused the analysis on treatment starting at the beginning of the infection. 

The following theorem determines the effectiveness of antivirals in terms of $\eta_p^c(t_{tr})$ and $\eta_{\beta}^c(t_{tr})$, by considering them separately. In the next subsection, the combined effect will be studied.

\begin{theorem}\label{theo:mainresult}
	Consider system \eqref{eq:SysOrigAcutAntiv}, constrained by the positive set $\X$. Consider that the virus spreads in the host for $t\ge 0$, which implies that: $U(0) > \Uc$, $I(0) = 0$ and $V(0) > 0$. Consider also that, at time $t_{tr}$, with $\check t_V < t_{tr} < \hat t_V $ (when the virus is increasing), an antiviral treatment is initiated such that $\eta_p(t)$ or $\eta_{\beta}(t)$ jumps from $0$ to $\eta_p \in [0,~1)$ or $\eta_{\beta} \in [0,~1)$, respectively (as stated in \eqref{eq:etas}). Then, 
	\begin{enumerate}
		\item if the inhibition effect is such that $\eta_p > \eta_p^c(t_{tr})$ (or $\eta_{\beta} > \eta_{\beta}^c(t_{tr})$) then the maximum of $V(t)$ occurs at $t_{tr}$, which is smaller than $\hat t_V$ by hypothesis. In this case \textbf{the antiviral treatment is effective}.
		\item for two antiviral treatment, 1 and 2, with inhibition effects $\eta_p^1,\eta_p^2 $ such that $\eta_p^c(t_{tr}) < \eta_p^1 < \eta_p^2 $ (or two inhibition effects $\eta_{\beta}^1,\eta_{\beta}^2 $ such that $\eta_{\beta}^c(t_{tr}) < \eta_{\beta}^1 < \eta_{\beta}^2 $), there exist a time
		$t^*>t_{tr}$ (large enough) such that $V_2(t) < V_1(t)$, for all $t\in (t_{tr},t^*]$, being $V_1(t)$ and $V_2(t)$ the virus corresponding to treatments $\eta_p^1$ and $\eta_p^2$, 
		respectively. In this case \textbf{both treatments are effective, and treatment 2 is more efficient than treatment 1}.
		\item if the inhibition effect is such that $\eta_p<\eta_p^c(t_{tr})$ (or $\eta_{\beta}<\eta_{\beta}^c(t_{tr})$), then the virus reaches a maximum at 
		a time $\hat t_{V,tr} >t_{tr}$. Furthermore, there exists a time $t^{e}=t^{e}(\eta_p,t_{tr})$ - denoted as early treatment time - such that $\hat t_{V,tr} > \hat t_{V}$ for $t_{tr} \in (\check t_V, t^{e}]$. In this case \textbf{the antiviral treatment is ineffective}.
		\item for two antiviral treatment, 1 and 2, with inhibition effects $\eta_p^1,\eta_p^2 $ such that $ \eta_p^1 < \eta_p^2 < \eta_p^c(t_{tr})$ with $t_{tr} \in (\check t_V, t^{e}(\eta_p,t_{tr})]$ (or for two inhibition effects $\eta_{\beta}^1,\eta_{\beta}^2 $ such that $ \eta_{\beta}^1 < \eta_{\beta}^2 < \eta_{\beta}^c(t_{tr})$), it is 
		$\hat t_{V,tr}^2 > \hat t_{V,tr}^1 >\hat t_V$, being $\hat t_{V,tr}^1$ and $\hat t_{V,tr}^2$ the virus maximum time corresponding 
		to treatments $\eta_p^1$ and $\eta_p^2$, respectively. In this case \textbf{both treatments are ineffective, but treatment 1 is more efficient than treatment 2}, which is a rather counter-intuitive fact.
		%In this case it is said that \textbf{the antiviral treatment is ineffective}.  
	\end{enumerate}
\end{theorem}
\begin{spf}
	For the sake of simplicity and clarity the proof is based on the approximation of system \eqref{eq:SysOrigAcut} described in Remark \ref{cor:Vmax}. In such a case, system \eqref{eq:SysOrigAcut} is approximated by $\dot U(t)=-\beta U(t) V(t)$, $\dot V(t) \approx (\Rn(t)-1) \delta V(t)$ and $I(t) \approx \frac{c}{p}V(t)$. Only antiviral treatments affecting parameter $p$ is considered, since the results follows from conditions on $\Rn(t)$ and parameters $\beta$ and $p$ affect $\Rn$ in the same way.	
\begin{enumerate}
		\item Since $\check t_V < t_{tr} < \hat t_V$ (by hypothesis) and $\hat t_V \approx t_c$, then by Theorem \ref{theo:key0}, $\Rn(t)>1$ for $t<t_{tr}$. Given that $\eta_p > \eta_p^c(t_{tr})$, then by \eqref{eq:Rtreat} $\Rn(t)$ jumps to a value smaller than one at $t_{tr}$, \textit{i.e.}, $\Rn(t_{tr}^-)>1$ and $\Rn(t_{tr})<1$ (where $\Rn(t_{tr}^-):=\lim_{t\to t_{tr}^-}\Rn(t)$). Given that $U(t)$ is strictly decreasing, $\Rn(t)<1$ for $t>t_{tr}$. From equation $\dot V(t) \approx (\Rn(t)-1) \delta V(t)$, it follows that $\dot V(t)<0$ for all $t>t_{tr}$. So, given that is assumed that $V(t)$ is increasing when the treatment is initiated, its maximum occurs at $t_{tr}$.
		\item Let us consider that for all $t\ge t_{tr}$, $V_i(t)$ and $U_i(t)$ are the virus and the susceptible cells for treatment $i=1,2$, \textit{i.e.}, the treatment under the inhibition effect $\eta_p^i$, and $\Rn_i(t) = \frac{U_i(t)(1-\eta_p^i) p\beta}{c \delta}$ with $U_1(t_{tr})=U_2(t_{tr})$ and $V_1(t_{tr})=V_2(t_{tr})$. 
		By hypothesis $\eta_p^c(t_{tr}) < \eta_p^1 < \eta_p^2$, then $\Rn_2(t_{tr})<\Rn_1(t_{tr})$ and $\Rn_i(t)<1$ for all $t\ge t_{tr}$ and $i=1,2$. 
		Since $U_i(t)$ is a continuous function for $i=1,2$, there is a positive time $t^*>t_{tr}$ such that $\Rn_2(t)<\Rn_1(t)$ for all $t\in [t_{tr},t^*]$.
		From equation $\dot V_i(t) \approx (\Rn_i(t)-1) \delta V_i(t)$ it follows that $\dot V_2(t)<\dot V_1(t)<0$ for all $t\in (t_{tr},t^*]$ (note that $\dot V$ is not defined on $t_{tr}$). As it is shown in Section \ref{sec:simres_1}, for real data patients, $t^*$ is large enough such that $V(t^*) \approx 0$.
		\item Since $\check t_V < t_{tr} < \hat t_V$ (by hypothesis) and $\hat t_V \approx t_c$, then by Theorem \ref{theo:key0}, $\Rn(t)>1$ for $t < t_{tr}$. Given that $\eta_p < \eta_p^c(t_{tr})$, then by \eqref{eq:Rtreat} $\Rn(t)$ jumps, at $t_{tr}$, from $\Rn(t_{tr}^-) = \frac{U(t_{tr}) p \beta}{c \delta}$ to the smaller value $\Rn(t_{tr}) = \frac{U(t_{tr}) p(1-\eta_p) \beta}{c \delta}$, which is still greater than one; \textit{i.e.}, $\Rn(t_{tr}^-) >\Rn(t_{tr})>1$. From the fact that $\Rn(t_{tr})>1$, it follows that $\dot V(t_{tr})>0$, and given that $\Rn(t)$ is decreasing for $t>t_{tr}$, there exists a time $\hat t_{V,tr}>t_{tr}$ such that $\Rn(\hat t_{V,tr})=1$, in which case, it is 
		\begin{eqnarray}\label{eq:theo1}
		\dot V(\hat t_{V,tr}) &\approx& (\Rn(\hat t_{V,tr})-1) \delta V(\hat t_{V,tr})=0,\\
		\ddot{V}(\hat t_{V,tr})&\approx&-\frac{p \beta^2}{c}U(\hat t_{V,tr})V(\hat t_{V,tr}) + (\Rn(\hat t_{V,tr})-1)^2 \delta^2 V(\hat t_{V,tr})<0.
		\end{eqnarray}
	This means that $V(t)$ reaches a maximum at time $\hat t_{V,tr} > t_{tr}$. Now, we need to prove that there exist a time $t^{e}(\eta_p,t_{tr})$ such that for $t_{tr} \in (\check t_V, t^{e})$ it is $\hat t_{V,tr} > \hat t_{V}$. According to Lemma \ref{lem:Rt}, there exists a time $t^{e} = t^{e}(\eta_p,t_{tr})$, smaller than $\hat t_V$, such that for a treatment time $t_{tr} \in (\check t_V, t^{e}]$, $\hat t_{V,tr}$ is a decreasing function of $\Rn(t_{tr})$, for $\Rn(t_{tr})>1$. This implies that smaller values of $\Rn(t_{tr})$ (greater values of $\eta_p$) correspond to larger times at which $\Rn(t)$, $t>t_{tr}$, reaches $1$,\textit{ i.e}., smaller values of $\Rn(t_{tr})$, correspond to larger values of $\hat t_{V,tr}$. So, from the fact that  $\Rn(t_{tr}^-) >\Rn(t_{tr})$ (being $\Rn(t_{tr}^-)$ the value of $\Rn$ at $t_{tr}$ if no treatment is applied), it follows that $\hat t_{V,tr} > \hat t_{V} > t_{tr}$, for $t_{tr} \in (\check t_V, t^{e}]$. As it is shown in Sections \ref{sec:simres_1} and \ref{sec:simres_2}, for real data patients, (the maximal) $t^e$ is close to $\hat t_V$.
		\item Since $\eta_p^1 < \eta_p^2 < \eta_p^c(t_{tr})$, by hypothesis, then $1<\Rn_2(t_{tr})<\Rn_1(t_{tr})$, being $\Rn_1(t_{tr}) = \frac{U(t_{tr}) p(1-\eta_p^1) \beta}{c \delta}$ and $\Rn_2(t_{tr}) = \frac{U(t_{tr}) p(1-\eta_p^2) \beta}{c \delta}$. Therefore, by following the same steps of the previous item, it follows that $\hat t_{V,tr}^2 > \hat t_{V,tr}^1 >\hat t_V$, being $\hat t_{V,tr}^1$ and $\hat t_{V,tr}^2$ the maximum times of the virus corresponding to $\eta_p^1$ and $\eta_p^2$, respectively. $\square$
	\end{enumerate}
\end{spf}

\begin{rem}
According to Theorem \ref{theo:key0}, $\Rn(t_{tr})$ should be compared to $1+\alpha(t_{tr})$, not to $1$, to obtain $\eta_p^c(t_{tr})$ and $\eta_{\beta}^c(t_{tr})$. However, as stated in Remark \ref{rem:alphaneg}, the correction factor $\alpha(t_{tr})$ is significantly smaller than $1$ for the considered period of time (when $t_{tr}$ is smaller than $\hat t_V$), so it was disregarded for the sake of simplicity. The effect of $\alpha(t_{tr})$ may be significant only for $t_{tr}$ close to $\hat t_V$ and, in such case, the viral dynamic can be slightly different from the one prescribed by Theorem \ref{theo:mainresult} (and Lemma \ref{lem:Rt}).
\end{rem}

\begin{rem}
	Item (iii) of Theorem \ref{theo:mainresult} establishes just the existence of $t^{e}>\check t_V$ such that for treatments starting at $t_{tr} \in (\check t_V,t^{e}]$, the new virus peak time is larger than the one corresponding to the untreated case, \textit{i.e.}, $\hat t_{V,tr} > \hat t_V$.  
	However, it should be noted that for parameters coming from real patient data, (the maximal) $t^{e}$ is indeed close to $\hat t_V$. This means that the time period where treatments can be ineffective is in most of the cases similar to $ (\check t_V,\hat t_V)$, \textit{i.e.}, similar to the time period where the virus is growing. Figures in Section \ref{sec:simres} confirm this fact.
\end{rem}

A main consequence of item (\textit{iii}) of Theorem \ref{theo:mainresult} is that early treatments, if no effective, are largely more dangerous than late ones.
Another critical point to be remarked is that when an early treatment is not strong enough to avoid the virus spreadability right after $t_{tr}$, the greater the antiviral effectiveness $\eta_p$ (or $\eta_{\beta}$) is, the more is the time the virus remains in the host, since the maximum time $\hat t_{V,tr}$ is delayed, as established in item (iv). As a result, antivirals could be detrimental as treated patients would need to be isolated for larger periods of time than untreated ones. 

Note that even when virus peak time can be delayed for some treatments, the virus peak will be always smaller than the one without any treatment. Furthermore, the fraction of dead cells at the end of the infection ($D$) will be always greater if no antiviral is administrated. Since $D=1-U_\infty/U(t_{tr})$, if we consider initial conditions at treatment time in equation \ref{eq:Uinf}, and being $\Uc^{tr} = U(t_{tr})/R(t_{tr})$, it results 
\begin{eqnarray}\label{eq:D}
    D=1 + \frac{W(-\Rn(t_{tr})e^{-(\Rn(t_{tr})+\K(t_{tr}))})}{\Rn(t_{tr})},
\end{eqnarray}
where $\K(t_{tr})=\frac{\beta}{c} (\frac{p(1-\eta_p)}{\delta} I(t_{tr}) - V(t_{tr}))=\Rn(t_{tr})I(t_{tr})-\frac{\beta}{c}V(t_{tr})$ and $W(\cdot)$ is the Lambert function ($W(\cdot)=W_p(\cdot)$, see \cite{abuin2020char}). Consequently, if treatment is started early, such that $\Rn(t_{tr}) \gg \K(t_{tr})$,  then the fraction of dead cells can be approximated as $D=1 + \frac{W(-\Rn(t_{tr})e^{-\Rn(t_{tr})})}{\Rn(t_{tr})}$, which is equal to 0 if $R(t_{tr})\leq 1$. Note that, by definition, $W(-\Rn(t_{tr})e^{-\Rn(t_{tr})})=-\Rn(t_{tr})$ for $R(t_{tr})\leq 1$, so $D=1+(-\Rn(t_{tr})/\Rn(t_{tr}))=0$. On the other hand, if $\Rn(t_{tr})>1$ then $W(-\Rn(t_{tr})e^{-\Rn(t_{tr})}) \in (-1,0]$ is monotonically increasing with $\Rn(t_{tr})$ (see Figure 1, \cite{abuin2020char}). Hence, the fraction of dead cells at the end of infection is a monotonically increasing function of $\Rn(t_{tr})$, for $\Rn(t_{tr})>1$. Moreover, as the treatment time is delayed, $\K(t_{tr})$ grows, being $D$ a monotonically increasing function of $\Rn(t_{tr})$, even for $\Rn(t_{tr}) \leq 1$.

%%%%%%%%%%%%%%%%%%%%%%%%%%%%%%%%%%%%%%%%%%%%%%%%%%%%%
\subsection{Antiviral effectiveness considering the combined effect on $\eta_\beta$ and $\eta_p$}\label{sec:andcharact}
%%%%%%%%%%%%%%%%%%%%%%%%%%%%%%%%%%%%%%%%%%%%%%%%%%%%%

\begin{figure}[H]
	\centering
	\includegraphics[width=0.7\textwidth]{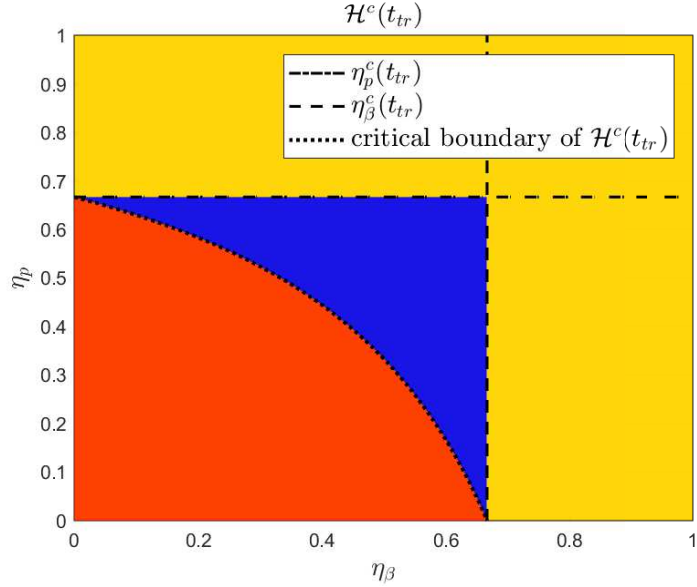}
	\caption{$\mathcal{H}^c(t_{tr})$ is given by the yellow and the blue regions, considering the following system parameters: $\beta=0.5$, $\delta=0.2$, $p=2$ and $c=5$, with $U(t_{tr})=3$. The critical boundary of $\mathcal{H}^c(t_{tr})$ represents the critical pairs of $\eta_\beta$ and $\eta_p$.}
	\label{fig:Hc}
\end{figure}

Theorem~\ref{theo:mainresult} describes the behavior of the virus under the effect of an inhibitor reducing the infection rate $\beta$ or the replication rate $p$. 
If both effects ($\eta_\beta$ and $\eta_p$) are simultaneously included in the model, it can be computed a region - in the $\eta_\beta,\eta_p$ space - for which condition $\Rn(t_{tr}) < 1$ is fulfilled. This way, instead of independent critical values $\eta_\beta^c(t_{tr})$ and $\eta_p^c(t_{tr})$ corresponding to each parameter, there is an entire set of critical combinations that makes a treatment effective. This set depends on treatment time $t_{tr}$ and is placed on the boundary of the effective set
\begin{equation}\label{eq:CriticalRegion}
\mathcal H^c(t_{tr}):=\{ (\eta_{\beta},\eta_p)\in [0,1]\times[0,1]: \eta_p > 1 - \frac{c\delta}{U(t_{tr})\beta p (1-\eta_{\beta})} \}.
\end{equation}
Clearly, every pair $(\eta_\beta,\eta_p)\in\mathcal{H}^c(t_{tr})$ fulfills condition $\Rn(t_{tr})<1$ and, reciprocally, every pair $(\eta_\beta,\eta_p)\notin\mathcal{H}^c(t_{tr})$ fulfills condition $\Rn(t_{tr})>1$. Figure \ref{fig:Hc} shows a plot of set $\mathcal{H}^c(t_{tr})$ in the plane $(\eta_{\beta} \times \eta_p)$, corresponding to parameters $\beta=0.5$, $\delta=0.2$, $p=2$ and $c=5$, with $U(t_{tr})=3$.  

As it can be seen, for every inhibition effect pair $(\eta_\beta,\eta_p)\not\in\mathcal{H}^c(t_{tr})$ (\textit{i.e.} a point inside the red region in Figure \ref{fig:Hc}) the antiviral treatment is ineffective. On the other hand, as it is shown in Theorem \ref{theo:mainresult}, for every pair $(\eta_\beta,\eta_p)$ such that $\eta_\beta>\eta_\beta^c(t_{tr})$ or $\eta_p>\eta_p^c(t_{tr})$ (the yellow region in Figure \ref{fig:Hc}) the antiviral treatment is effective. Finally, for pairs in the region inside $\mathcal{H}^c(t_{tr})$ with $\eta_\beta<\eta_\beta^c(t_{tr})$ and $\eta_p<\eta_p^c(t_{tr})$ (which does not fit conditions in Theorem \ref{theo:mainresult} and it is represented by the blue region
in Figure \ref{fig:Hc}) the antiviral treatment is still effective.
%%%%%%%%%%%%%%%%%%%%%%%%%%%%%%%%%%%%%%%%%%%%%%%%%%%%%%%%%%%%%%%%%%%%%%%%%%%%%%%%%%

%%%%%%%%%%%%%%%%%%%%%%%%%%%%%%%%%%%%%%%%%%%%%%%%
\section{Within-Host Modeling of COVID-19}\label{sec:pac}
%%%%%%%%%%%%%%%%%%%%%%%%%%%%%%%%%%%%%%

In this section, the parameters of model \eqref{eq:SysOrigAcut} are estimated using viral load data of 9 RT-PCR COVID-19 positive patients, labeled as A,B,C,D,E,F,G,H and I, reported by Woelfel et. al. \cite{wolfel2020virological}. We follow a similar procedure as in Vargas et. al. \cite{vargas2020host}. Since the viral load is measured in Log10 scales, the model fitting was fulfilled by minimizing the root mean squared of logarithmic error (RMSLE), denoted as: 
\begin{equation}\label{eq:Param_Est_Cost}
RMSLE = \sqrt{\frac{1}{n}\sum_{i=1}^{n}(log(y_i)-log(\bar{y}_i))^2},
\end{equation}
where $n$ is the number of measurements, $y_i$ the model predictive output, and $\bar{y}_i$ the experimental measurement. Since the minimization of \ref{eq:Param_Est_Cost} implies a nonlinear optimization problem, with highly dependence on initial conditions, the Differential Evolution (DE) algorithm \cite{storn1995differential,hernandez2019modeling} is employed as a global optimization algorithm, which has shown to be robust to initial guesses of parameters \cite{torres2016comparative}. 

\begin{figure}[hbt!]
	\hspace{-1.5cm}
	\includegraphics[trim = 35mm 0mm 35mm 0mm, clip, width=1.2\textwidth]{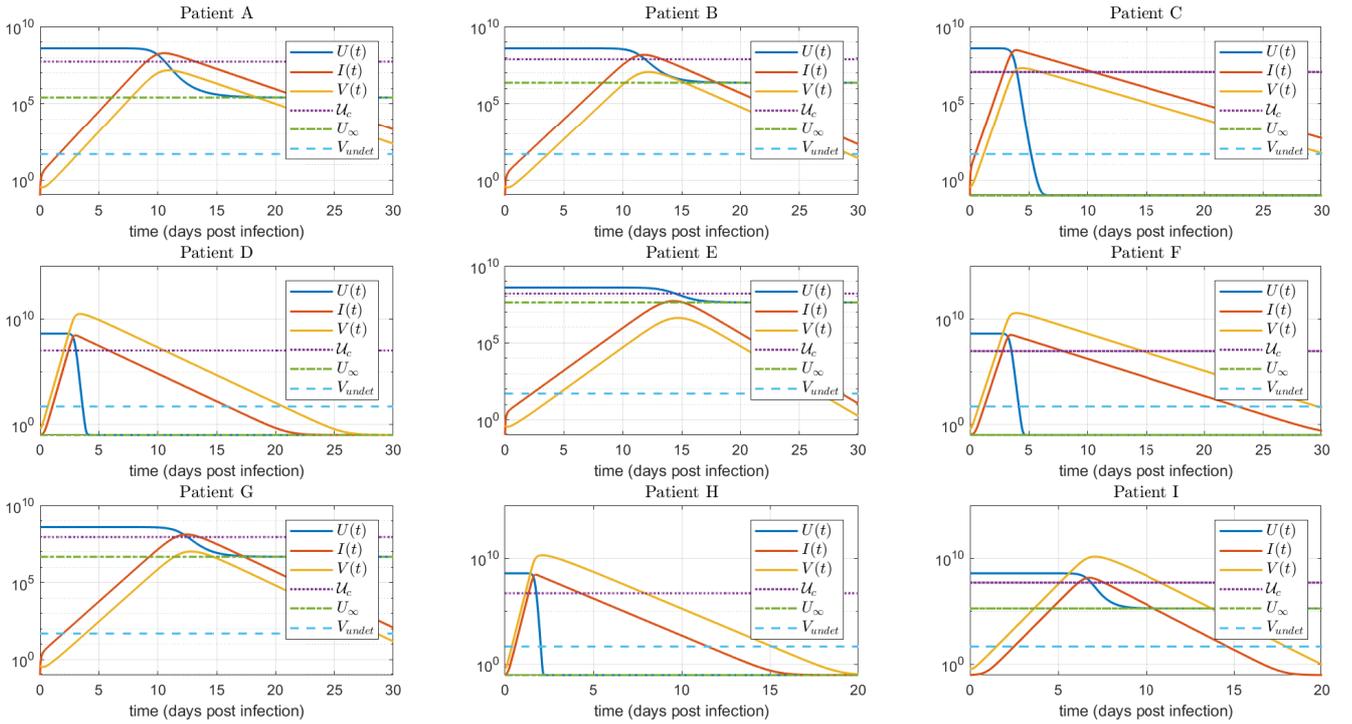}
	\caption{\small{$U$,$I$,$V$ time evolution for the untreated case. As stated in Theorem \ref{theo:key0}, $\check t_V < \hat t_I < t_c < \hat t_V$, where $\check t_V$ and $\hat t_V$ are the times at which $V(t)$ reaches a local minimum and a local maximum, respectively, $\hat t_I$ is the time at which $I(t)$ reaches a local maximum, and $t_c$ is the time at which $U(t)$ reaches $\Uc$. Furthermore, as stated in Remark \ref{cor:Vmax}, $\check t_V \approx 0$ and $\hat t_I \approx t_c \approx \hat t_V$. $V_{undet} = 50$ (Copies/mL) stands for the undetectable level of the virus load.}}
	\label{fig:allpatOL}
\end{figure}

Even though it is still debatable which compartments SARS-CoV-2 can infect, there is a common agreement that the viral shedding take places mainly in the respiratory epithelial cells (due to the high expression of ACE2)  with direct viral toxicity of the infected cells \cite{gupta2020extrapulmonary}. Therefore, following previous works of mathematical modelling for influenza infection in humans, the value of $U_0$ is taken as about $4 \times 10^8$ cells for all patients \cite{baccam2006kinetics}. Furthermore, $I_0$ is assumed to be $0$ and
$V_0$ was estimated (using a regression model, since the viral at the day of infection was not provided)  \cite{vargas2020host} to be about $0.31$ Copies/mL. Moreover, in order to avoid identifiability problems related with the fact that only viral titters are used to fit the model, the viral clearance parameter (c) was set in $2.4$ 1/day, which is in accordance with previous estimates for influenza and HIV \cite{hernandez2013modeling,baccam2006kinetics}. 
The parameters and initial condition (at the time of infection $t=0$) of each patient are collected in Table \ref{tab:param}. Since the viral load was measured after the onset of the symptoms, an incubation period of 7 days post infection (dpo) was assumed for the time of infection, according to \cite{vargas2020host}.
\begin{table}[H]
	\begin{center}		
		\caption{Target cell-limited model parameter values for COVID-19 patients \cite{vargas2020host}}
		\begin{tabular}{| c | c | c | c |c |}% c|c|c|}
			\hline
%			\multicolumn{5}{ |c| }{Amount of drug (mg)} \\ \hline
	    Patient       &  $\beta$ &  $\delta$ &  $p$  & $c$ \\\hline%&  $U_0$ &  $I_0$ & $V_0$  \\ \hline
		A  &  $1.35\times10^{-7}$  &  0.61 &  0.2  &  2.4 \\%&	$1.0\times10^{7}$      &  0  &  5.01 \\%
		B  &  $1.26\times10^{-7}$  & 0.81 &  0.2  &  2.4\\%&	$1.0\times10^{7}$      &   0  &  0.31 \\
		C  &  $5.24\times10^{-7}$  & 0.51 & 0.2  &  2.4\\%&	$1.0\times10^{7}$      &   0  &  0.31	\\
		D  &  $7.92\times10^{-10}$  & 1.21 & 361.6  &  2.4\\%&	$1.0\times10^{7}$      &   0  &  0.31	\\
		E  &  $1.51\times10^{-7}$  &  2.01 &  0.2  &  2.4\\%&	$1.0\times10^{7}$      &   0  &  0.31	\\
		F  &  $5.74\times10^{-10}$  & 0.81 & 382  &  2.4\\%&	$1.0\times10^{7}$      &   0  &  0.31	\\
		G  &  $1.23\times10^{-7}$  &  0.91 & 0.2  &  2.4\\%&	$1.0\times10^{7}$      &   0  &  0.31	\\
	    H &	  $2.62\times10^{-9}$   & 1.61 & 278.2  &  2.4\\%&	$1.0\times10^{7}$      &   0  &  0.31	\\
        I &	  $3.08\times10^{-10}$  &  2.01 & 299   & 2.4\\\hline%&	$1.0\times10^{7}$      &   0  &  0.31	\\ \hline
		\end{tabular}
		\label{tab:param}
	\end{center}
\end{table}
The time evolution of $U$, $I$ and $V$ is shown in Figure \ref{fig:allpatOL}, for each patient. As it can be seen, the plot confirms the results in Theorem \ref{theo:key0} concerning the minimum and maximum times of $V$, the maximum time of $I$ and the time when $U$ reaches $\Uc$ (\textit{i.e.}, when $\Rn(t)$ reaches $1$).
Figure \ref{fig:allpatOL} also confirms the approximation made in Remark \ref{cor:Vmax}, by showing that virus and infected cell peaks occur approximately at the same time at which $U$ reaches $\Uc$.
All these time values, together with the values of $\Uc$ (the critical target cell value), $U_{\infty}$ (the final value of $U$, when the virus is died out),
$\RN$ and the virus maximum are shown in Table \ref{tab:const}, for the nine simulated patients.
\begin{table}[H]
	\begin{center}		
		\caption{Parameter Characterization of COVID-19 patients. The times are given in days post infection (dpi).}
		\begin{tabular}{| c | c | c | c | c | c | c | c | c | c | c | c |}
			\hline
			%			\multicolumn{5}{ |c| }{Amount of drug (mg)} \\ \hline
			Patient & $\Uc~(\times 10^8)$ &$U_{\infty}$ & $\RN$ & $\check t_V$ & $\hat t_I$ &  $t_c$ & $\hat t_V$ & $V_{max}$ \\ \hline
			A		&	0.54        & 	0.26 $\times 10^6$  & 7.35  &   0.22       & 10.53    &  10.61 &  10.95     &  1.50 $\times 10^7$ \\
			B		&	0.77        & 	0.23 $\times 10^7$  & 5.17  &   0.24       & 11.81    &  11.90 &  12.20     &  1.18 $\times 10^7$  \\
			C		&	0.12        & 	0.00      & 34.22 &   0.06       & 3.94     &  3.99  &  4.42      &  2.10 $\times 10^7$ \\
			D		&   0.10        &   0.00      & 39.47 &   0.02       & 3.00     &  3.05  &  3.39      &  2.90 $\times 10^{10}$\\
			E		&   1.60        &	4.27 $\times 10^7$  & 2.50  &   0.23       & 14.31    &  14.47 & 14.69      &  0.42 $\times 10^7$  \\
			F		&	0.09        &	0.00      & 45.12 &   0.03       & 3.48     &  3.53  &  3.93      &  3.55 $\times 10^{10}$ \\
			G		&	0.90        &	0.47 $\times 10^7$  & 4.50  &   0.25       & 12.44    &  12.53 &  12.83     &  1.04 $\times 10^7$ \\
			H		&	0.05        &	0.00      & 75.57 &   0.00       & 1.73     &  1.77  &  2.11      &  2.02 $\times 10^{10}$ \\
			I		&	0.52        & 	0.19 $\times 10^6$  & 7.64  &   0.07       & 6.77     &  6.86  &  7.09      &  1.50 $\times 10^{10}$  \\  \hline
		\end{tabular}
		\label{tab:const}
	\end{center}
\end{table}

\begin{figure}[hbt!]
	\hspace{-1.5cm}
	\includegraphics[trim = 35mm 0mm 30mm 0mm, clip, width=1.15\textwidth]{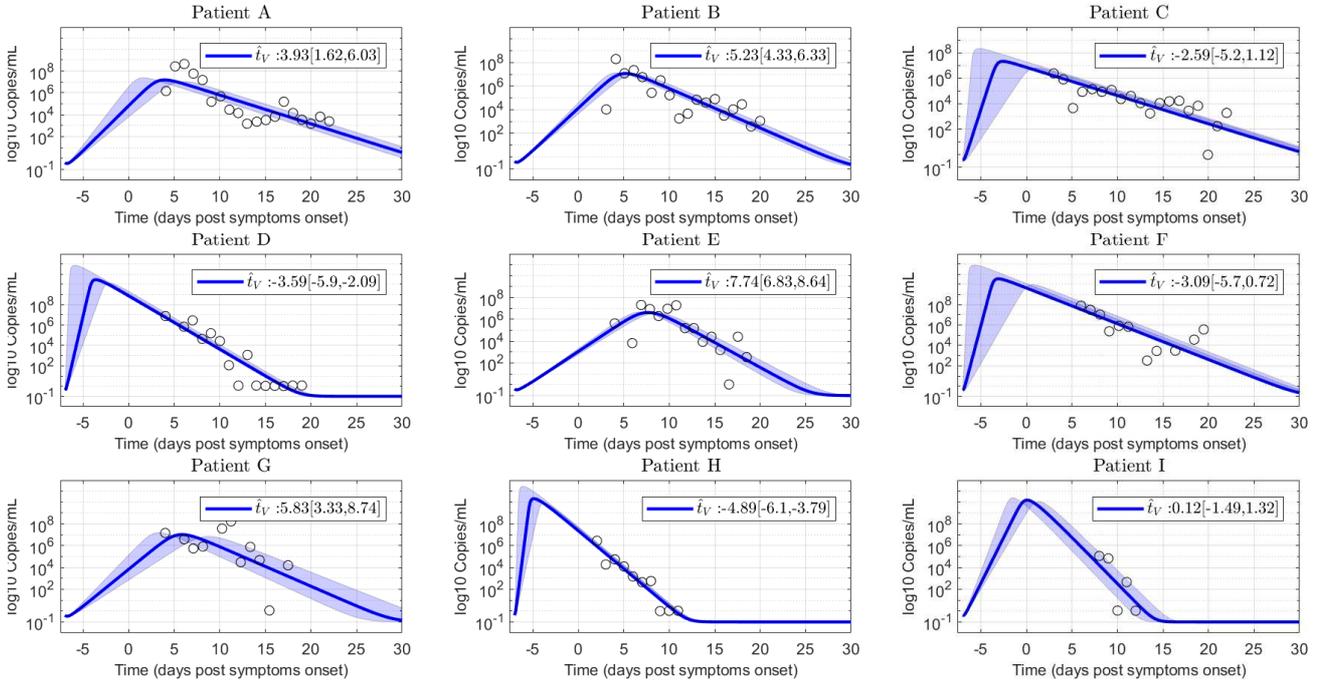}
	\caption{\small{Uncertainty analysis on viral load evolution for $p$-parameter $95\%$ confidence interval. Note the uncertainty in viral load peak. Empty dots are data of COVID-19 patients.}}
	\label{fig:SA_Casep}
\end{figure}

\begin{rem}
    In comparison with other coronavirus diseases, such as MERS and SARS, where the virus load peak takes place after the onset of symptoms ($7-10$ days post infection) \cite{he2020temporal}, for SARS-CoV-2 infection it is not clear the temporal interval where the viral load reaches a peak. A recent study linking epidemiological and viral load data, suggests that the viral load peak occurs during the day of symptom onset \cite{he2020temporal}. However, observation of viral load in infected macaques \cite{williamson2020clinical,munster2020respiratory} denotes that the viral peak from nose and throat swabs happens during the 1-3 days post infection. Therefore, since the target cell model fitting was conducted using SARS-CoV-2 viral load measured after the onset of the symptoms \cite{wolfel2020virological}, the estimated time-to-peak ($\hat t_V$) in Table \ref{tab:const} is subject to practical identifiability problems, which as  was indicated in Theorem \ref{theo:mainresult} and will be shown in Section \ref{sec:simres}, is a main parameter to evaluate antiviral effectiveness. Due to this reason, an uncertainty analysis computing the likelihood-based confidence intervals \cite{raue2009structural} for each parameter was done. In Figure \ref{fig:SA_Casep} such analysis is shown for the $95\%$ confidence interval of parameter p, where it can be seen the high degree of uncertainty in the estimated $\hat t_{V}$
    for patients A, C, F and G (approximately $5$ days). 
\end{rem}
%%%%%%%%%%%%%%%%%%%%%%%%%%%%%%%%%%%%%%%%%%%%%%%%%%%%%%%%%%%%%%%%%%%%%%%%%%%%%%%%%%

%%%%%%%%%%%%%%%%%%%%%%%%%%%%%%%%%%%%%%%%%%%%%%%%%%%%%%%%%%%%%%%%%%%%%%%%%%%%%%%%%%
\section{Simulation Results} \label{sec:simres}
%%%%%%%%%%%%%%%%%%%%%%%%%%%%%%%%%

To evaluate the results concerning the antiviral treatment effectiveness, several simulations involving the nine patients introduced in Section \ref{sec:pac} were performed. First, the virus spreading interval - the time between the estimated day of infection ($7~ dpo$) and the time-to-peak of viral load ($\hat t_V$) - will be considered in Section \ref{sec:simres_1}, to assess the results of Theorem \ref{theo:mainresult}. Due to the reported variability on the estimated time-to-peak (Figure \ref{fig:SA_Casep}), 
%happening before the symptoms onset for 4/9 patients (labeled as patient C, D, F and H), 
we decided to initialize the treatments taking into account the relative time with respect to the estimated $\hat t_V$, instead of a fixed time for all the population. Furthermore, the (maximal) early treatment time ($t^e$) of each patient was computed numerically, being of the order of $0.77 ~[0.74,0.79] \hat t_V$.

Then, in Section \ref{sec:simres_2} the case when the antiviral therapy is started before and after the untreated time-to-peak $\hat t_V$ is simulated, in order to analyze the sub-potent/potent drug efficacy as treatment time $t_{tr}$ is delayed. Finally, in Section \ref{sec:simres_3} the synergistic effects of antiviral therapies blocking the viral replication $p$ and the host cell infection $\beta$ are studied, taking into account the combined drug effect analysis made in Section \ref{sec:andcharact}. 

To numerically assess the viral kinetics evolution the following infected-related metrics are employed: 
i) the difference of viral loads with and without treatment at time-to-peak, $\Delta V=V(\hat t_V)-V(\hat t_{V,tr})$, which is a measure of the viral reduction at time-to-peak respect to the untreated case \cite{gonccalves2020timing}; 
ii) the duration of infection, $DI$, defined as the time spent by the viral titer curve over a detection limit of 100 copies/ml, which is a measure of the viral shedding interval, and 
iii) the time-to-peak $\hat t_V$, which is an indicator of the viral replication rate. 
For the sake of simplicity, unless otherwise stated, the approximation of system \eqref{eq:SysOrigAcut} presented in Remark \ref{cor:Vmax} is used for simulations. %As stated in Section \ref{sec:pd}, $\alpha(t_{tr})$ will be ignored, since it is smaller than $10^{-2}$ for most of the patients and treatment times. 

%%%%%%%%%%%%%%%%%%%%%%%%%%%%%%%%%
\subsection{Treatment initiated at different times, during viral spreading interval ($t_{tr}\in (\check t_V,\hat t_V)$)} \label{sec:simres_1}
%%%%%%%%%%%%%%%%%%%%%%%%%%%%%%%%%
\vspace{0.2cm}
%%%%%%%%%%%%%%%%%%%%%%%%%%%%%%%%%
\subsubsection*{Scenario 1: Treatment is initiated at viral load detection level ($t_{tr}=t_{DL}$)}
%%%%%%%%%%%%%%%%%%%%%%%%%%%%%%%%%
\vspace{0.2cm}

\begin{figure}[hbt!]
	\hspace{-1.5cm}
	\includegraphics[trim = 32mm 0mm 32mm 0mm, clip, width=1.2\textwidth]{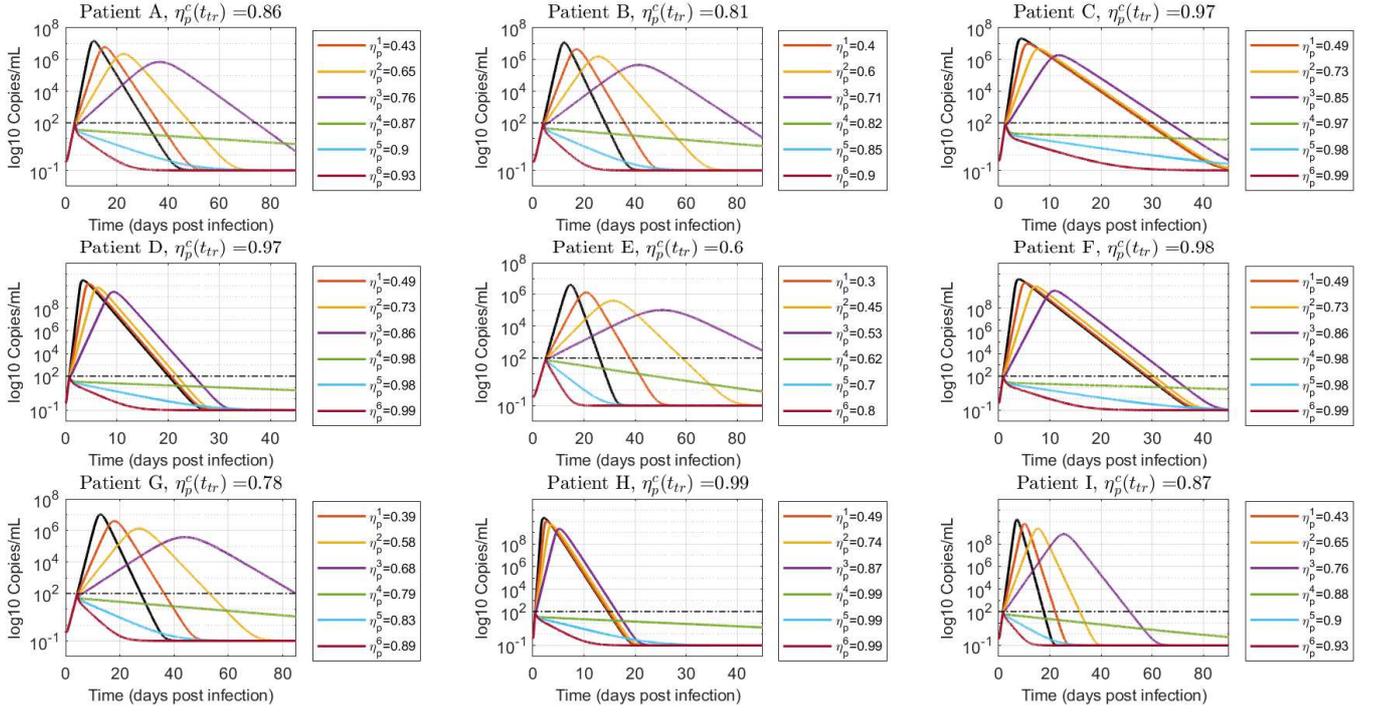}
	\caption{\small{Free virus behavior when treatment is started at viral load detection level ($t_{tr}=t_{DL}$). Values of $\eta_p$ smaller ($\eta_p^1,\eta_p^2$), approximately equal ($\eta_p^3,\eta_p^4$) and greater ($\eta_p^5,\eta_p^6$) than $\eta_p^c$ are simulated to demonstrate the results in Theorem \ref{theo:mainresult}. The black line denotes the untreated case ($\eta_p=0$).}}
	\label{fig:allpatCL_DL}
\end{figure}

Figure \ref{fig:allpatCL_DL} shows the simulated virus load evolution corresponding to each of the nine patients, when the antiviral therapy is started at the first positive PCR test day ($t_{tr}=t_{DL}$, being $t_{DL}$ the time at which the virus reaches the value $DL=100$ copies/mL), which is about $0.3\hat t_V \approx 0.40 t^e$. Consequently, all the patients fall in the case $\check t_V < t_{tr} < t^{e}<\hat t_V$ which means that we are under the hypothesis of Theorem \ref{theo:mainresult} and, furthermore, the treatment time belongs to the early treatment time interval. To properly asses the antiviral replication inhibition effect, increasing values of $\eta_p$ were used. First we started with small values of $\eta_p$, fulfilling the condition $\eta_p < \eta_p^c(t_{tr})$ (item iii in Theorem \ref{theo:mainresult}, ineffective treatment). Then $\eta_p$ was set to a value equal to $\eta_p^c(t_{tr})$, to reinforce the counter-intuitive fact that larger inhibition effects fulfilling $\eta_p < \eta_p^c(t_{tr})$ produce larger virus peak times (item iv in Theorem \ref{theo:mainresult}). Finally, $\eta_p$ was increased to two values larger than $\eta_p^c(t_{tr})$ (items i and ii in Theorem \ref{theo:mainresult}), to show that higher inhibition effects produce faster eradication of the virus $V$. 

As expected, effective treatments produce an instantaneous decline of the viral load, with a monotonically decreasing viral shedding interval as the antiviral efficacy is incremented. On the other hand, ineffective treatments cause a delay in time to peak, significantly increasing the duration of viral shedding as $\eta_p$ is augmented from $0$ to $\eta_p^c$. Even when the viral load peak is a monotonically decreasing function of $\eta_p$ (by following similar steps than the ones in Lemma \ref{lem:Rt} it can be shown that $\partial V(\hat t_{V,tr})/\partial \Rn= (\delta/\beta)(1-\frac{1}{\Rn})$, for $\Rn:=\Rn(t_{tr})\geq1$, with $\Rn(t_{tr})=(1-\eta_p)p\beta/(c\delta)$), the patient will be PCR-positive for longer periods of time. This means that isolation and precautions measures with treated patients should be carefully considered, depending on the antiviral effectiveness.   

\begin{figure}[hbt!]
    \begin{center}
	\hspace{-1.5cm}
	\includegraphics[width=0.75\textwidth]{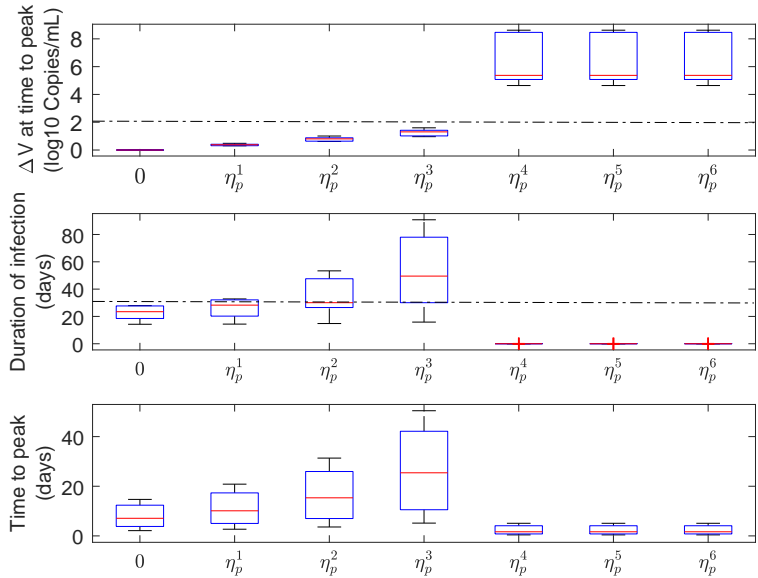}
	\caption{\small{Infection-related metrics as function of $\eta_p$ ($t_{tr}=t_{DL}$), for antiviral effectiveness assessment.}}
	\label{fig:Indexes_allpatCL_DL}
	\end{center}
\end{figure}

Figure \ref{fig:Indexes_allpatCL_DL} shows a box-plot of the infected-related metrics for the antiviral effectiveness assessment. For an effective antiviral therapy ($\eta_p>\eta_p^c(t_{tr})$), the difference of viral loads at time-to-peak ($\Delta V$) is above the 2 logs threshold \cite{gonccalves2020timing} and the duration of infection ($DI$) is below the 30 days limit (according \cite{zhou2020duration} the viral shedding interval for untreated COVID-19 patients is in order of 30 days) which is in accordance with an effective viral clearance strategy. On the other hand, for ineffective antiviral therapy ($\eta_p<\eta_p^c(t_{tr})$), even tough $\Delta V$ is monotonically increasing with $\eta_p$, the duration of infection is increased to a sub-potent drug efficacy. This means that the delay of the virus peak associated to ineffective treatments is significant in terms of the infected-related metrics.

%

%%%%%%%%%%%%%%%%%%%%%%%%%%%%%%%%%
\subsubsection*{Scenario 2: treatment initiated in the course of viral spreading}
%%%%%%%%%%%%%%%%%%%%%%%%%%%%%%%%%

\begin{figure}[hbt!]
	\hspace{-1.5cm}
	\includegraphics[trim = 32mm 0mm 32mm 0mm, clip, width=1.2\textwidth]{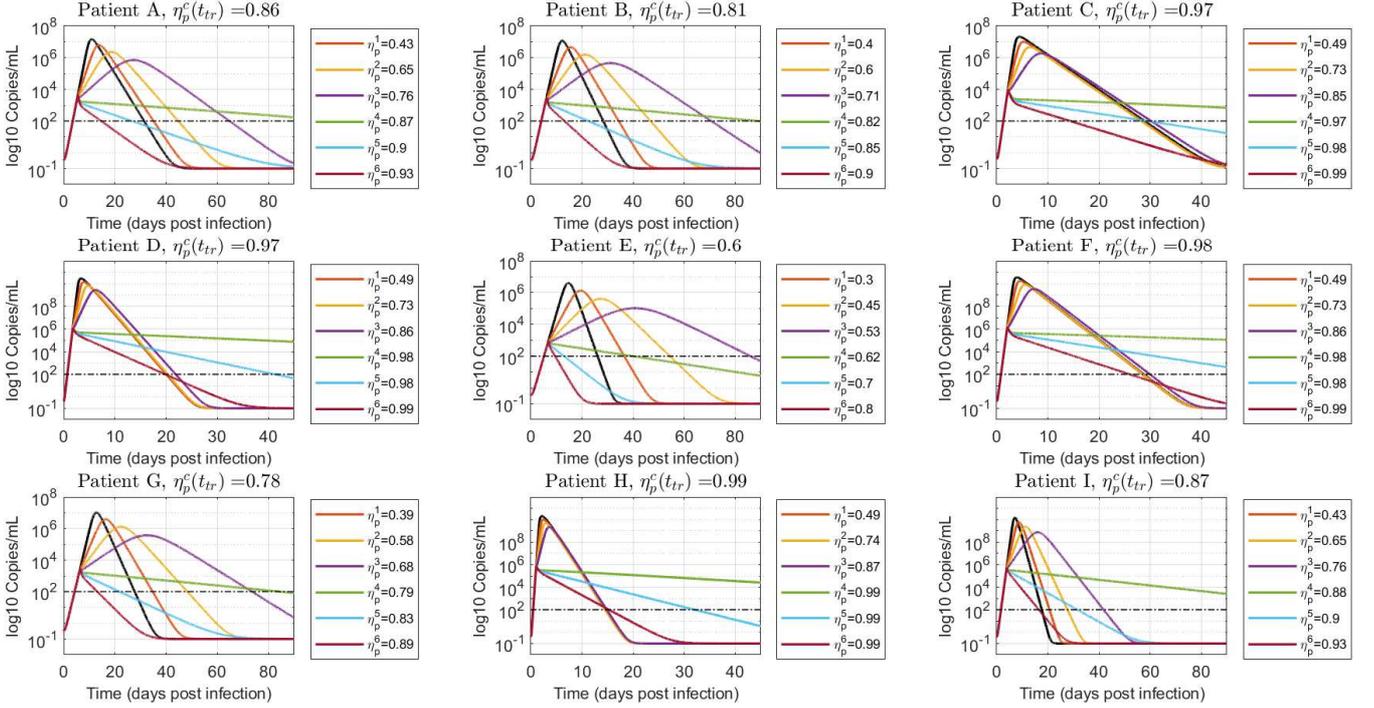}
	\caption{\small{Viral load time evolution with treatment initial time given by $t_{tr}=0.7t^{e}$. Values of $\eta_p$ smaller, approximately equal and greater than $\eta_p^c$ are simulated to demonstrate the results in Theorem \ref{theo:mainresult}.The black line denotes the untreated case ($\eta_p=0$).}}
	\label{fig:allpatCL_07te}
\end{figure}

\begin{figure}[hbt!]
	\hspace{-1.5cm}
	\includegraphics[trim = 32mm 0mm 32mm 0mm, clip, width=1.2\textwidth]{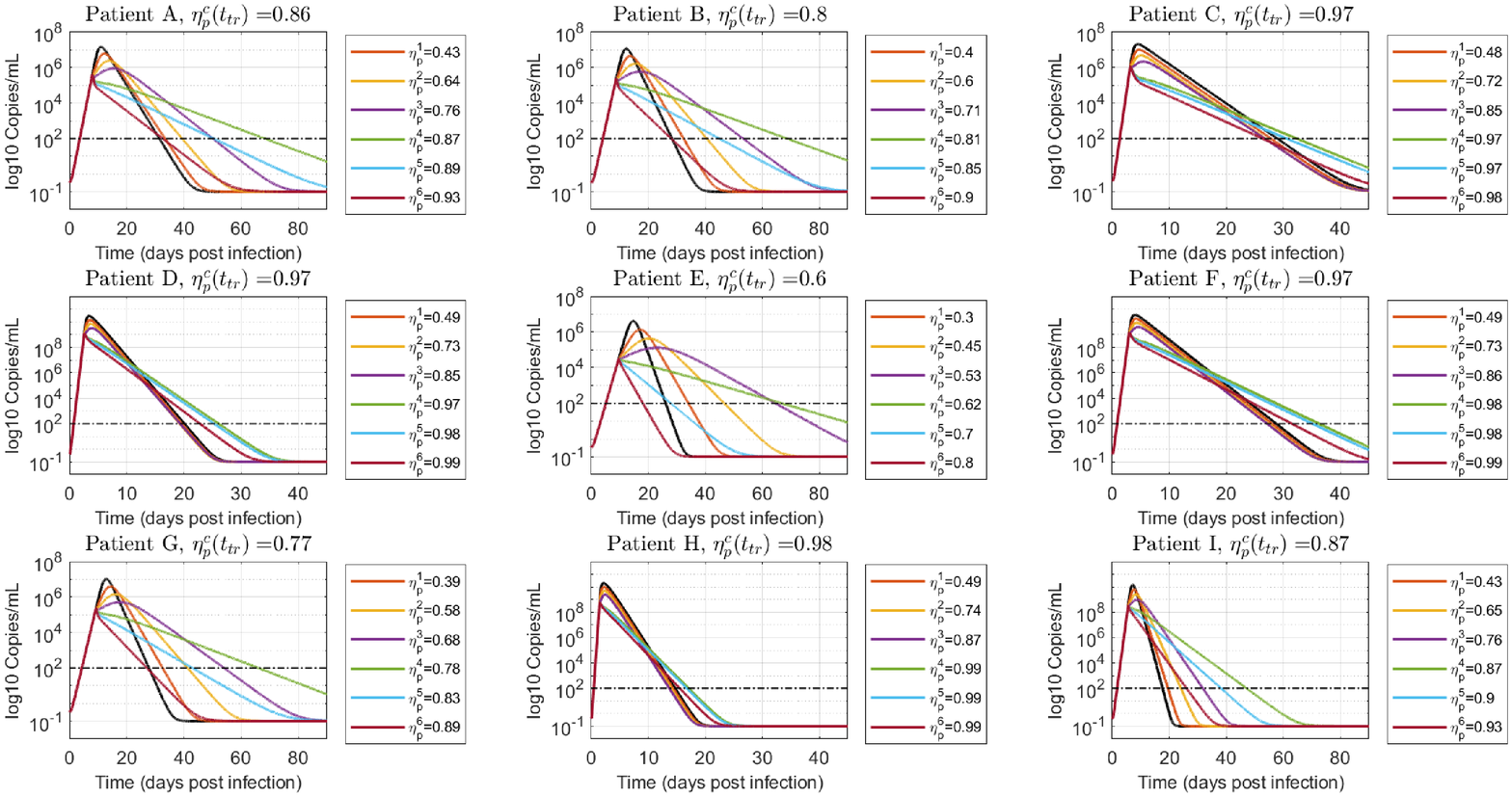}
	\caption{\small{Viral load time evolution with treatment initial time given by $t_{tr}=t^{e}$. Values of $\eta_p$ smaller, approximately equal and greater than $\eta_p^c$ are simulated to demonstrate the results in Theorem \ref{theo:mainresult}. The black line denotes the untreated case ($\eta_p=0$).}}
	\label{fig:allpatCL_te}
\end{figure}

\begin{figure}[hbt!]
    \begin{center}
	\hspace{-1.5cm}
	\includegraphics[width=0.75\textwidth]{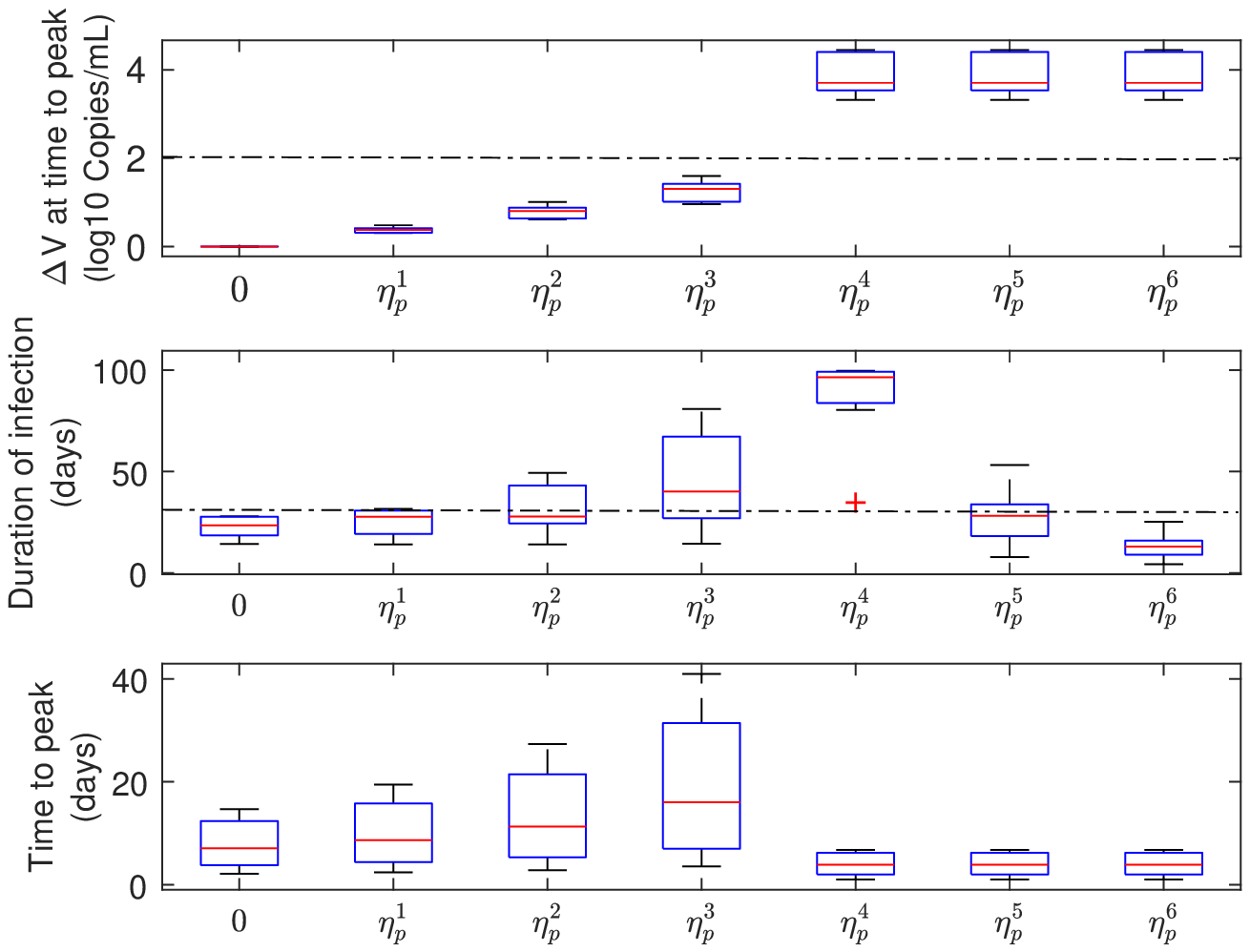}
	\caption{\small{Infection-related metrics as function of $\eta_p$ ($t_{tr}=0.7t^{e}$), for antiviral effectiveness assessment.}}
	\label{fig:Indexes_allpatCL_07te}
	\end{center}
\end{figure}

\begin{figure}[hbt!]
    \begin{center}
	\hspace{-1.5cm}
	\includegraphics[width=0.75\textwidth]{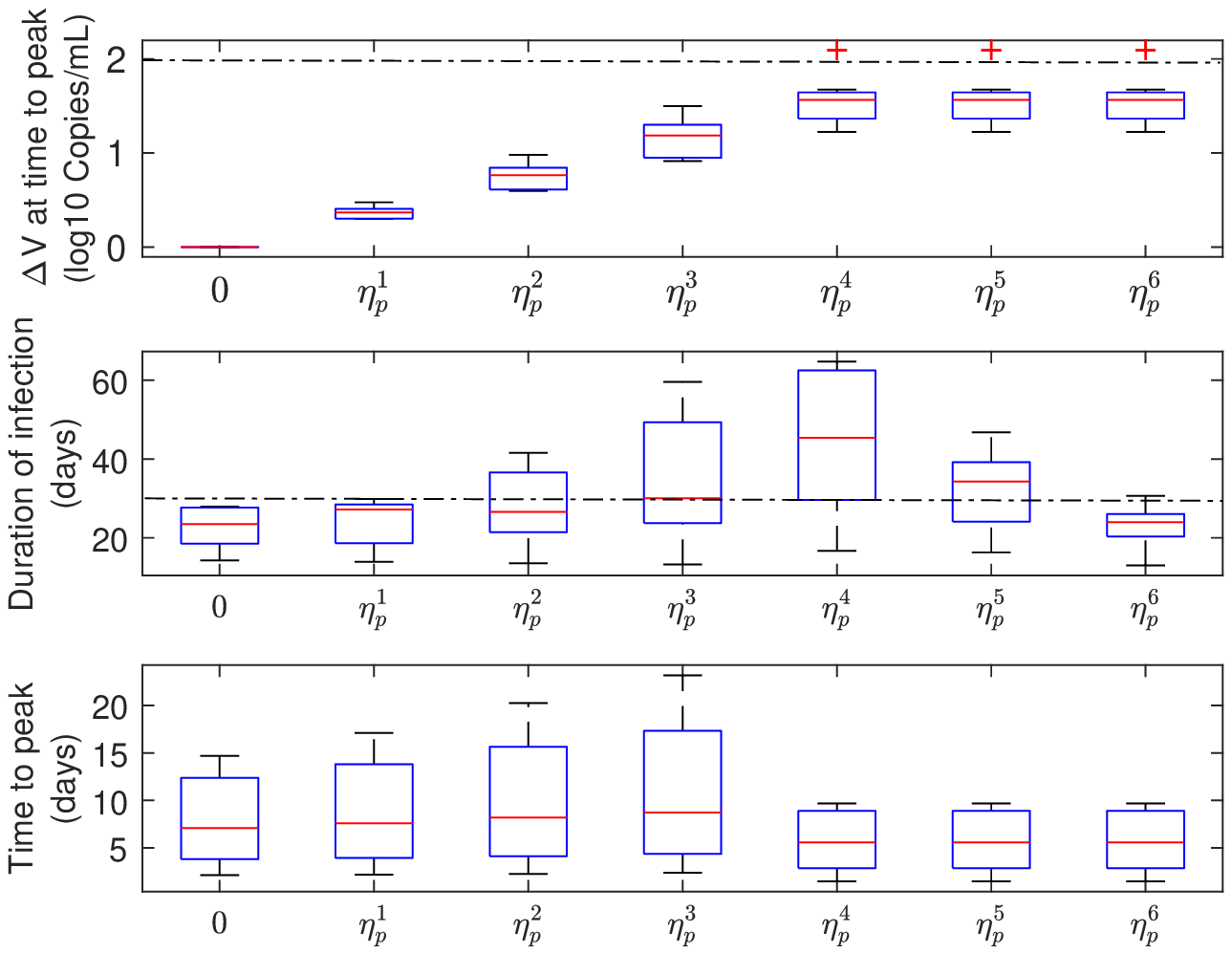}
	\caption{\small{Infection-related metrics as function of $\eta_p$ ($t_{tr}=t^{e}$), for antiviral effectiveness assessment.}}
	\label{fig:Indexes_allpatCL_te}
	\end{center}
\end{figure}

Figures \ref{fig:allpatCL_07te} and \ref{fig:allpatCL_te} show the virus load evolution when the antiviral therapy is started in the course of viral spreading, at $t_{tr}\approx 0.7t^{e}$ and $t_{tr}\approx t^{e}$, respectively. The viral load over the time has the same qualitative behavior than in the previous case. However, the duration of the infection is longer as the inhibition effect approaches $\eta_p^c(t_{tr})$ from below and from above (\textit{i.e.}: $\eta_p^3$ and $\eta_p^4$ cases, violet and green lines in Figures \ref{fig:allpatCL_07te} and \ref{fig:allpatCL_te}). Indeed, for values of drug effectiveness in the vicinity of $\eta_p^c(t_{tr})$,
%$\Rn(t_{tr})\approx 1$, and consequently, since $\dot V(t)\approx (\Rn(t_{tr})-1) \delta V(t)$ for $t>t_{tr}$, 
the viral duration interval is augmented - even for effective therapies - since $V$ has already reached a relatively high value at the treatment time. This behavior is confirmed in box-plots \ref{fig:Indexes_allpatCL_07te} and \ref{fig:Indexes_allpatCL_te}, where a sudden increase of the duration of infection happens when the antiviral efficacy is near the critical value. However, as the treatment time is delayed this behavior is mitigated (see box-plot \ref{fig:Indexes_allpatCL_te}) due to the natural increment of $V(t_{tr})$ (recall that $\dot V(t)\approx (\Rn(t_{tr})-1) \delta V(t)$, for $t>t_{tr}$, so larger values of $V(t_{tr})$ produce larger values of $|\dot V(t)|$). It is important to remark that, even when an increased duration of the viral shedding is reported for the effective treatment case in the boundary of critical drug effectiveness ($\eta_p\approx\eta_p^c$), the viral load peak occurs previous to $\hat t_V$, in accordance with items i and ii of Theorem \ref{theo:mainresult}. 

%%%%%%%%%%%%%%%%%%%%%%%%%%%%%%%%%
\subsection{Treatment initiated at different times, with the same effectiveness} \label{sec:simres_2}
%%%%%%%%%%%%%%%%%%%%%%%%%%%%%%%%%

\begin{figure}[hbt!]
	\hspace{-1.5cm}
	\includegraphics[trim = 32mm 0mm 32mm 0mm, clip, width=1.2\textwidth]{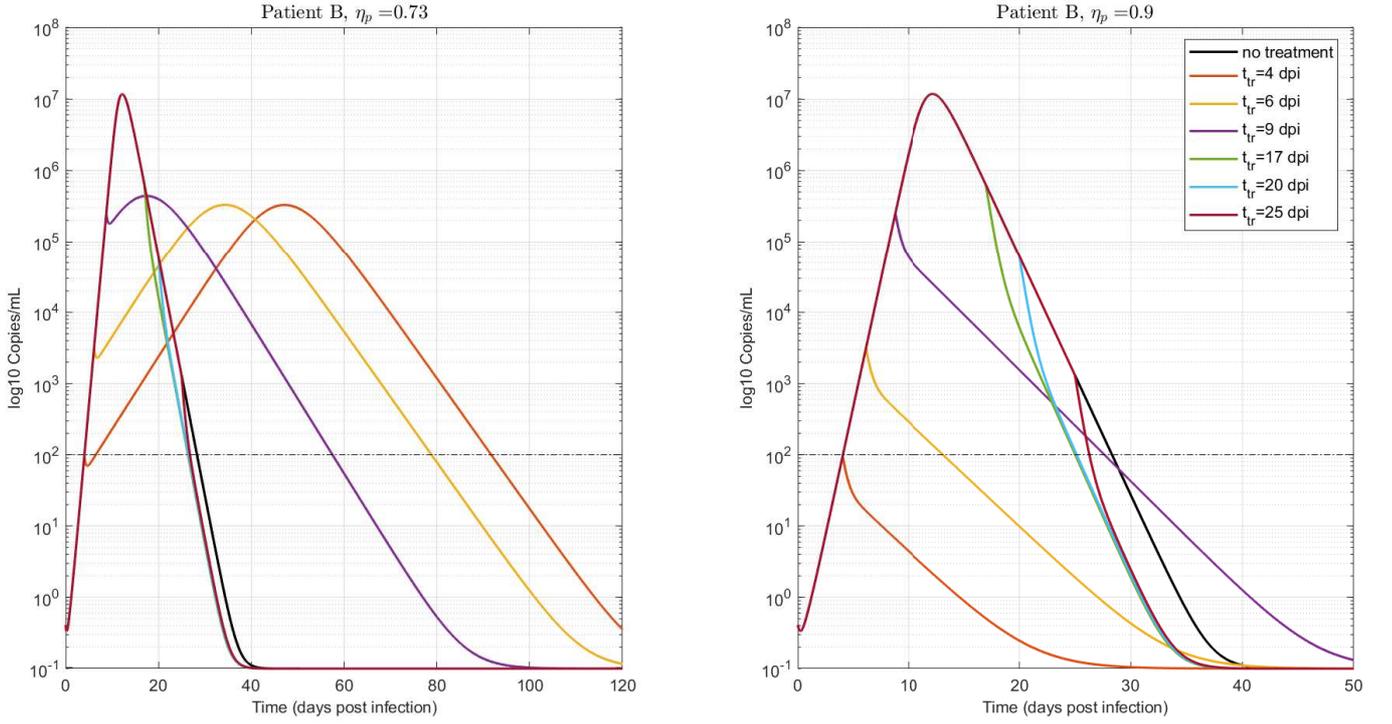}
	\caption{\small{Virus time evolution for different treatment times, $t_{tr}=4,6,8,17,20,25$ (dpi). Two fixed values of 	$\eta_p$ were used, smaller and bigger than $\eta_p^c(t_{tr})$: $\eta_p=0.73$ (left) and $\eta_p=0.9$ (right), respectively. $\eta_p^c(t_{tr})\approx 0.81$ for $t_{tr}<t^e$. Patient B. }}
	\label{fig:patBCL_tr}
\end{figure}

\begin{figure}[hbt!]
	\hspace{-1.5cm}
	\includegraphics[trim = 32mm 0mm 32mm 0mm, clip, width=1.2\textwidth]{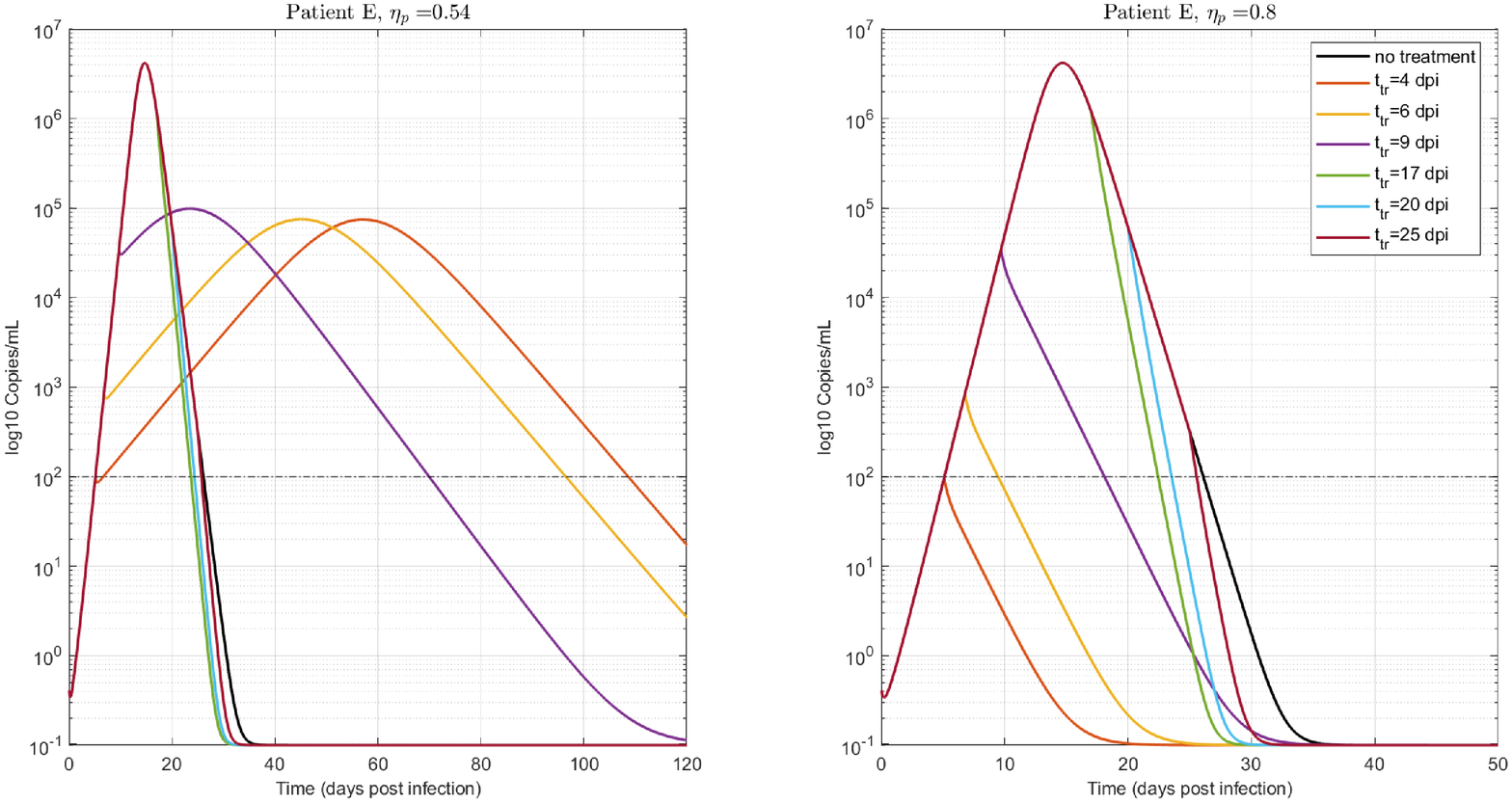}
	\caption{\small{Virus time evolution for different treatment times, $t_{tr}=4,6,8,17,20,25$ (dpi). Two fixed values of $\eta_p$ were used, smaller and bigger than $\eta_p^c(t_{tr})$: $\eta_p=0.54$ (left) and $\eta_p=0.8$ (right), respective.  $\eta_p^c(t_{tr})\approx 0.6$ for $t_{tr}<t^e$. Patient E.}}
	\label{fig:patECL_tr}
\end{figure}

In order to analyze the viral kinetics when the antiviral therapy is initiated at any time in the whole infection period (\textit{i.e.}, before and after $\hat t_V$), we studied the temporal dependence of infection-related metrics corresponding to fixed sup-potent/potent drug effectiveness. Only two representative patients (labeled as B and E) were considered, based on the uncertainty analysis made in Section \ref{sec:pac}. The treatment initial times were: $t_{tr1} = 4$, $t_{tr2} = 6 $, $t_{tr3} = 9$, $t_{tr4} = 17$, $t_{tr5} = 20$ and $t_{tr6} = 25$ dpi. Figures \ref{fig:patBCL_tr} and \ref{fig:patECL_tr} (left) show that if an ineffective treatment is initiated ($\eta_{pB}=0.73$ for patient B and $\eta_{pE}=0.54$ for patient E, with $\eta_{pB}^c(t_{tr})\approx0.81$ and $\eta_{pE}^c(t_{tr})\approx 0.6$ for $t_{tr}<t^e$), the time-to-peak ($\hat{t}_{V,tr}$) and the duration of infection ($DI$) decrease as the treatment time is delayed. This can be explained by the fact that, for the same effectiveness, $\hat{t}_{V,tr}$ is reduced as $V(t_{tr})$ is augmented. Equations \eqref{eq:lemm6} and \eqref{eq:lemm6_1} in Lemma \ref{lem:Rt} show that for a fixed $\Rn\geq1$, $k_1(\Rn)\geq 0$ and $k_2(\Rn)\leq 0$ are increasing functions of $V(t_{tr})$ (\textit{i.e.}:  $k_1(\Rn)=0$ and $k_2(\Rn)=-1/(\delta(\Rn-1))$ for $V(t_{tr})\rightarrow 0$ while $k_1(\Rn)=1/((1-\Rn)^2 +\Rn)$ and $k_2(\Rn)=-1/(\delta\Rn)$ for $V(t_{tr})\rightarrow \delta/\beta$, being $\delta/\beta<V(\hat t_V)$). Therefore, since $\hat t_{V,tr} =k_2 \ln (k_1) + t_{tr}$ (according to equation \eqref{eq:lemm4}), it can be inferred that $\hat{t}_{V,tr}$ is monotonically decreasing with $V(t_{tr})$. Moreover, comparing the viral load at the time-to-peak, it can be deduced that $V(\hat t_{V,tr2})>V(\hat t_{V,tr1})$ for an ineffective therapy started at two treatment times $t_{tr2}>t_{tr1}$ during the beginning of the infection (assuming that $U(t_{tr1})\approx U(t_{tr2})$). 

By following similar steps as in Lemma \ref{lem:Rt}, the viral load at the time-to-peak can be written as: $V(\hat t_{V,tr})=V(t_{tr})+(p/c)U(t_{tr})-(\delta/\beta) (\ln (\Rn_{tr}) +1)$ and, therefore, $V(\hat t_{V,tr2})-V(\hat t_{V,tr1}) = V(t_{tr2})-V(t_{tr1})>0$. However, from Figures \ref{fig:patBCL_tr} and \ref{fig:patECL_tr} (right) it can be seen that $V(\hat t_{V,tr})$ is of the order of $10^6$ and $10^5$ copies/ml for patients B and E, respectively, independently of the treatment initiation time. This can be explained by the fact that the difference of viral load at $\hat t_{V,tr}$, for therapies started at different times, depends mainly on the deviation between the viral loads at treatment times, which normally are several order of magnitude below $V(\hat t_{V,tr})$. For example, for patient B, $V(\hat t_{V,tr1})=3.29\times 10^5$, $V(\hat t_{V,tr2})=3.30\times 10^5$ and $V(\hat t_{V,tr3})=4.35\times 10^5$ Copies/ml while $V(t_{tr1})=10^2$, $V(t_{tr2})=2.5 \times 10^3$ and $V(t_{tr3})=1.74 \times 10^5$ Copies/ml. Consequently, the increment on the viral peak is not significant as the treatment is delayed, although the duration of infection is considerably reduced (\textit{i.e.}: from $96$ to $55$ days for patient B and $95$ to $65$ days for patient E). Note that the target cell-model assumes that the viral clearance in the convalescent phase is proportional to the virus load concentration at time-to-peak, which is in the same order of magnitude for the three cases. As a result, and taking into account that the viral peak is reached early as the treatment is delayed, a quickly viral deletion is observed for the postponed case. As a consequence, if a fixed sub-potent drug is employed, delaying the treatment initiation time could reduce viral shedding interval, without significantly increasing the viral load at time-to-peak. It is important to remark that if $\eta_p$ is reduced even more ($\eta_p \ll \eta^c$), the viral at time-to-peak follows the same behavior, although, the viral shedding interval is not-considerably decreased, since $\hat t_{V,tr}\approx \hat t_{V}$ as $\eta \rightarrow 0$ ($t_{tr} \in (\check t_V,\hat t_V)$. Hence, although a non-significant improvement on the duration of infection is achieved if the treatment is delayed, the viral outcomes continues to be acceptable.

On the other hand, if the ineffective therapy is initiated after the viral load peak, the duration of infection is slightly reduced with respect to the previous case, since $R(t_{tr})<1$ at this time and, therefore, $\dot V=(\Rn(t_{tr})-1)\delta V(t)$ is a decreasing function of $V(t_{tr})$. In contrast to this, Figures \ref{fig:patBCL_tr} and \ref{fig:patECL_tr} (right) show that if an effective antiviral therapy is applied ($\eta_{pB}=0.9$ for patient B and $\eta_{pE}=0.8$ for patient E) the duration of the infection is still decreased, since
$\Rn_{1}(t_{tr4})<\Rn_{2}(t_{tr4})$ for both, patient B and patient E, being $\Rn_{i}$ the within-host basic reproduction number under a treatment with inhibition effect $\eta_p^i$, being $\eta_p^1=0.9$ and $\eta_p^2=0.73$ for patient B and $\eta_p^1=0.8$ and $\eta_p^2=0.54$ for patient E. 

%
%\todo{Ale. G. En las Figs, habria que agrandar las legends, y poner $\eta_p^c(tr1) =0.81$, etc. Si son todas iguales, sacarlas, y poner en Caption el unico valor}

%%%%%%%%%%%%%%%%%%%%%%%%%%%%%%%%%
\subsection{Treatment with antiviral inhibiting both viral infection ($\eta_{\beta}$) and viral replication rate ($\eta_{p}$) }
\label{sec:simres_3}
%%%%%%%%%%%%%%%%%%%%%%%%%%%%%%%%%

Finally, the virus behavior was considered when a combined antiviral therapy 
inhibiting both, the viral infection rate ($\beta$) and the 
viral replication rate ($p$), was applied. Since the critical drug efficacy for acute models is defined during the viral spreading interval (Section \ref{sec:pd}), the treatment was initiated at $t_{tr}=0.7 t^e$, which is a value fulfilling $t_{tr}\in(\check t_V, \hat t_V)$.

\begin{figure}[hbt!]
    \begin{center}
	\hspace{-1.65cm}
	\includegraphics[width=1.1\textwidth]{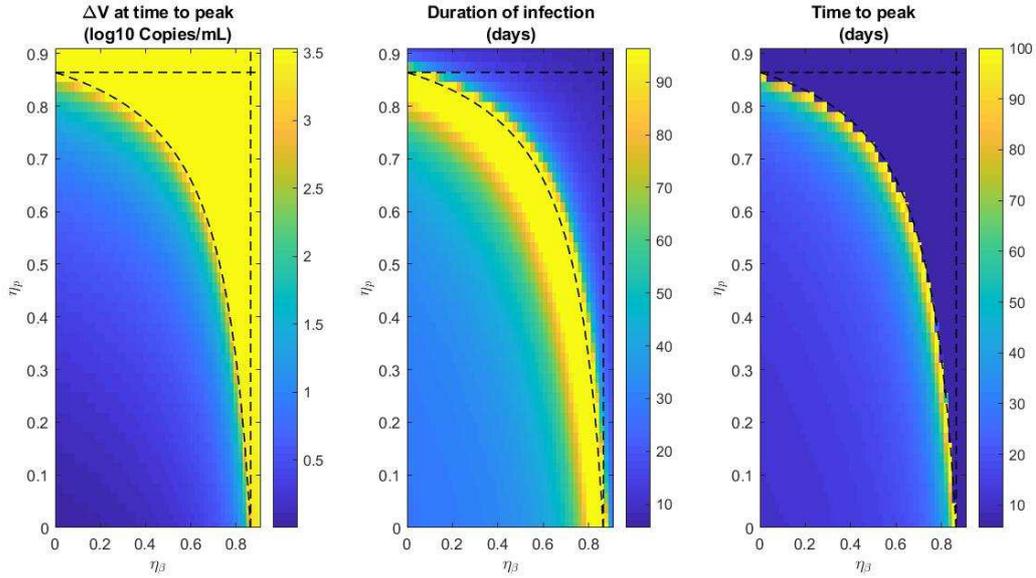}
	\caption{\small{Infection-related metrics for antiviral effectiveness assessment as function of $\eta_p$ and $\eta_{\beta}$ ($t_{tr}=0.7t^{e}$). Patient A}}
	\label{fig:Index_patACL_07te}
	\end{center}
\end{figure}

For the sake of clarity, only Patient A is considered and the 3 previous infection-related metrics are assessed as functions of the antiviral inhibition effects $\eta_p$ and $\eta_{\beta}$. Figure \ref{fig:Index_patACL_07te} shows that the antiviral is effective for every inhibition effect pair $(\eta_\beta,\eta_p)$ fulfilling $(\eta_\beta,\eta_p)\in\mathcal{H}^c(t_{tr})$ (\textit{i.e.}: yellow region of $\Delta V$), in accordance with the results in Section \ref{sec:andcharact}. Moreover, as for the single inhibition cases, the duration of infection and the time-to-peak increase as the combined drug efficacy pair $(\eta_\beta,\eta_p)$ tends to the boundary of $\mathcal{H}^c(t_{tr})$. The difference of viral load at time-to-peak $\Delta V$ remains constant for every $(\eta_\beta,\eta_p) \in \mathcal{H}^c(t_{tr})$ since $V(\hat t_{V,tr})=t_{tr}$ for effective treatments. In conclusion, the synergistic effect of drug effectiveness in combined therapies produces a reduction on the critical effectiveness with respect to single-therapies cases, with the same dynamical behavior over the boundary of $\mathcal{H}^c(t_{tr})$. 
%%%%%%%%%%%%%%%%%%%%%%%%%%%%%%%%%%%%%%%%%%%%%%%%%%%%%%%%%%%%%%%%%%%%%%%%%%%%%%%%%%
\section{Discussion}\label{sec:discuss}
%%%%%%%%%%%%%%%%%%%%%%%%%%%%%%%%%%%%%%%%%%%%%%%%%%%%%

While several vaccines have been developed to prevent COVID-19 disease, it is imperative to evaluate potential therapies against SARS-CoV-2 infection. Among them, antiviral treatments are promissory strategies to increment the viral clearance "within-host", decreasing disease severity. 
In this sense, the drug effectiveness concept is crucial to evaluate the drug effect threshold above which the viral load starts declining \cite{callaway2002hiv}. 
Although this critical value has been described for chronic infections (e.g. HIV \cite{callaway2002hiv,huang2003modeling}, HCV \cite{dahari2007modeling}) 
it has not been studied yet for COVID-19. 
The critical drug efficacy can be understood in terms of the antiviral potency to decline viral spreading, in spite of guiding the system to an uninfected equilibrium, as it was introduced in Definition \ref{def:anttreat}, Section \ref{sec:andcharact}. In addition, the critical inhibition effect of an antiviral depends on the treatment initiation time, and it has shown to be a monotonically decreasing function during the viral spreading interval (roughly speaking, during the time interval before the untreated virus peak).
%

% 7.2-Infected-related metrics when the antiviral therapy is started at different points of the viral spreading interval $t_{tr}\in (\check t_V,\hat t_V)$. 
Regarding the virus behavior for effective therapies ($\eta>\eta^c$), the viral load at time-to-peak is practically equal to $V(t_{tr})$ (\textit{i.e.}, a treatment is considered effective if the virus starts to decrease at the very moment the treatment is initiated, by modifying either $\eta_p$, $\eta_{\beta}$, or both), so $\Delta V$ is larger for earlier treatments, being greater than the 2 log threshold for effective antiviral therapies started at $t_{tr}<t^e$ ($t^e \approx 0.75 \hat t_V$). Nevertheless, for $\eta$ approaching $\eta^c$ (from above), the viral shedding interval is enlarged as the treatment is initiated during the infection growth. Therefore, if the therapeutic objective is to reduce both, the viral shedding interval and the viral peak, the efficacy level would need to overpass the critical value in a given quantity. From Figure \ref{fig:eta_c_tr} it can be seen that $\eta^c(t_{tr})>0.7$ for most of patients.

Simulations results suggest that for ineffective therapies, the virus would take longer to be cleared - rather than a monotonic decline, as in the previous case - reaching a peak and finally decreasing to zero. In spite of the fact that the viral load at time-to-peak decreases monotonically with $\eta$ ($\Delta V$ increases as the inhibition effect jumps from $0$ to $\eta >0$), the time-to-peak increases as $\eta$ approaches $\eta^c$ (from below), which results in a longer duration of infection, potentially requiring additional isolation measures for the treated patient.  

Although an effective reduction of the viral load peak and the duration of infection could be achieved with early treatments (treatments started before the untreated peak), minor effects are attained with late ones. Particularly, Figures \ref{fig:patBCL_tr} and \ref{fig:patECL_tr} corroborate a slight reduction on the viral shedding interval if a sub-potent drug therapy is employed after the time-to-peak. Consequently, taking into account that the viral load in COVID-19 patients presumably reaches the peak prior to the symptom onset, no further clinical improvements may be obtained if the therapy is started in this symptomatic phase
(\textit{i.e.}: notice that SARS-CoV-2 pathophysiology is characterized by a direct cytotoxic effect, endothelial cell damage, dysregulation of immune response, among others  \cite{gupta2020extrapulmonary}). Furthermore, considering the analysis presented in Theorem \ref{theo:mainresult} that is valid during viral spreading interval, we can argue that the non-significant antiviral effect reported recently in \cite{who2020repurposed} might be related to an early time-to-peak of COVID-19 patients. 

The effectiveness of combined treatments affecting both, $\eta_p$ and $\eta_{\beta}$, was also studied and an interdependent critical drug efficacy level was computed. Mathematically, the critical combination of values of $\eta_p$ and $\eta_\beta$ is placed on the boundary of the effective set $\mathcal{H}^c(t_{tr})$, which is a set in the $\eta_p,\eta_{\beta}$ space. In comparison with single treatment cases (represented by the horizontal and vertical dash-dot lines in Figure \ref{fig:Index_patACL_07te}), a reduced critical antiviral efficacy was reported for the combined case, denoting a reduction on the necessary drug effectiveness to reduce viral spreading. Moreover, the viral characteristic behavior in the vicinity of the critical drug efficacy (\textit{i.e.}: increase of time-to-peak and duration of infection) was preserved for the combined case, as the pair $(\eta_{\beta},\eta_{p})$ belongs to the critical boundary of $\mathcal{H}^c(t_{tr})$, which implies that the antiviral effectiveness characterization made in Theorem \ref{theo:mainresult} for single treatments could be extended to combined therapies.  

In sum, this work formalizes the existence of a critical drug efficacy for acute infection models, which could have implications in the extended viral shedding observed "in-silico" by \cite{goyal2020potency} (for SARS-CoV-2 infections) and by \cite{cao2017mechanisms} (for influenza) when the antiviral therapy is initiated early but with a sub-potent drug efficacy. Moreover, it was shown the importance of initialization the antiviral therapy early (before viral load peak) in order to achieve a significant reduction of $\Delta V$ and $DI$ \cite{gonccalves2020timing,wu2020early}. Although a time dependence was noticed for the critical drug efficacy ($\eta_p^c(t),\eta_{\beta}^c(t)$), since it is a decreasing function of $t$ (Figure \ref{fig:eta_c_tr}), its behavior does not compromise the antiviral success if an effective therapy is started later. 

The main clinical implications for acute infections, and, particularly, for the SARS-CoV-2 infection are: a) importance of viral load monitoring on probably infected COVID-19 patients (prophylactic use of antiviral therapy, previous to onset of symptoms, although adverse effects have been reported for potential antiviral drugs \cite{jomah2020clinical}, which could limited their prophylactic usage in risk patients), b) isolation of treated patient, c) possible explanation of the increase of duration of infection showed in immune compromised COVID-19 patients. \\ 
%%%%%%%%%%%%%%%%%%%%%%%%%%%%%%%%%%%%%%%%%%%%%%%%%%%%%%%%%%%%%%%%%%%%%%%%%%%%%%%%%%

%%%%%%%%%%%%%%%%%%%%%%%%%%%%%%%%%%%%%%%%%%%%%%%%%%%%%%%%%%%%%%%%%%%%%%%%%%%%%%%%%%
\newpage
\section{Appendix 1. Technical lemmas} \label{sec:app1}

%% LEMMA2
\begin{lem}\label{lem:Rt}
	Consider \eqref{eq:SysOrigAcutAntiv}, constrained by the positive set $\X$, at the beginning of the infection $t=0$, with
	$U(0) > \Uc$, $I(0) \geq 0$ and $V(0) > 0$. 
	Consider that at some time $t_{tr}>0$, $\eta_p(t)$ jumps from $0$ to $\eta_p > \eta_p^c$, being $\eta_p^c$ the critical value defined in \eqref{eq:eta_p}. Then, there exists $t^{e}(\eta_p,t_{tr})$ such that if $ t_{tr} \in (\check t_V, t^{e}]$, the virus peak time $\hat t_{V,tr}$ (considering the treatment effect) is a decreasing function of $\Rn(t_{tr})$.
\end{lem}

\begin{pf}
Consider the approximation \ref{eq:SysOrigAcutApp}, in Remark \ref{cor:Vmax}. Since $\eta_p(t)$ jumps from $0$ to $\eta_p < \eta_p^c$, at $t_{tr}$, then $\Rn(t)$ jumps from $\Rn(t_{tr}^-)$ to $\Rn(t_{tr})$, with $\Rn(t_{tr}) \in (\Rn(t_{tr}^-), 1)$ and $\Rn(t_{tr}^-)>1$. Then it is possible to approximate the explicit solutions for $U(t)$ and $V(t)$ (by approximating the $\ln(z)$ function in the time interval before the virus peak by $(z-1)(\frac{w_1}{z} + w_2)$, being $w_1$ and $w_2$ arbitrary constants, fulfilling $w_1+w_2=1$ \cite{kaushal2020estimating}. So, $U(t)$, for $t>t_{tr}$, can be written as
\begin{eqnarray}\label{eq:lemm1}
U(t) = \frac{(B + \sqrt{B^2-4 A C}) k_{t_{tr}}^N e^{-(\beta/U(t_{tr})) \sqrt{B^2-4 A C}(t-t_{tr})}-k_{t_{tr}}^D (B - \sqrt{B^2 -4 A C})}{2 A k_{t_{tr}}^D-2 A k_{t_{tr}}^N e^{-(\beta/U(t_{tr})) \sqrt{B^2-4 A C} (t-t_{tr})}}, 
\end{eqnarray}
where
\begin{eqnarray}\label{eq:lemm2}
k_{t_{tr}}^N &=& 2 A U(t_{tr}) + B - \sqrt{B^2 -4 A C}, \nonumber \\
k_{t_{tr}}^D &=& 2 A U(t_{tr}) + B + \sqrt{B^2 -4 A C}, \nonumber \\
A &=& -\frac{p(1-\eta_p)}{c} U(t_{tr}) + \frac{\delta}{\beta} w_2, \nonumber \\
B &=& U(t_{tr}) (\frac{\delta}{\beta} w_1-\frac{\delta}{\beta} w_2 + V(t_{tr}) + \frac{p(1-\eta_p)}{c}U(tr)), \nonumber \\
C &=& -\frac{\delta}{\beta} U(t_{tr})^2 w_1. \nonumber 
\end{eqnarray}
A reasonable approximation, however, can be obtained by selecting $w_1=0$ and $w_2=1$ (also $w_1=0.2$ and $w_1=0.8$ gives a good result, but clearly the former selection significantly simplifies the expressions). In such a case, $A$, $B$, $C$, $k_{t_{tr}}^N$ and $k_{t_{tr}}^D$ read: 
\begin{eqnarray*}\label{eq:lemm3}
A &\approx& -U(t_{tr})\frac{p(1-\eta_p)}{c} + \frac{\delta}{\beta} = (1-\Rn(t_{tr})) \frac{\delta}{\beta} <0, \nonumber \\
B &\approx& U(t_{tr}) (-\frac{\delta}{\beta} + V(t_{tr})+U(t_{tr})\frac{p(1-\eta_p)}{c})=-U(t_{tr}) A +U(t_{tr})V(t_{tr})>0, \nonumber \\
C &\approx& 0, \nonumber \\
k_{t_{tr}}^N &\approx& 2A U(t_{tr}), \nonumber \\
k_{t_{tr}}^D &\approx& 2 A U(t_{tr}) + 2B =  2 U(t_{tr})V(t_{tr}), \nonumber 
\end{eqnarray*}
being $\Rn(t_{tr}) = U(t_{tr}) \frac{\beta p(1-\eta_p)}{c \delta}$. This way, $U(t)$ can be simplified as
\begin{eqnarray}\label{eq:lemm1N}
U(t) &=& \frac{2B 2A U(t_{tr}) e^{-(\beta/U(t_{tr})) B (t-t_{tr})}}{2 A 2 U(t_{tr})V(t_{tr})-2 A 2A U(t_{tr}) e^{-(\beta/U(t_{tr})) B (t-t_{tr})}}, \nonumber \\
     &=& \frac{(U(t_{tr})V(t_{tr})-U(t_{tr})A) e^{-(\beta/U(t_{tr})) B (t-t_{tr})}}{V(t_{tr})-A e^{-(\beta/U(t_{tr})) B (t-t_{tr})}}
\end{eqnarray}

The time at which $U(t)$ reaches $\mathcal{U}_c^{tr}$ - \textit{i.e.}, the time at which $\Rn(t)$ reaches $1$ and $V(t)$ reaches its peak, denoted as $\hat t_{V,tr}$ (with $\hat t_{V,tr} \geq t_{tr}$) - can be explicitly computed as
\begin{eqnarray}\label{eq:lemm4}
\hat t_{V,tr} =k_2 \ln (k_1) + t_{tr},
\end{eqnarray}
where 
\begin{eqnarray*}\label{eq:lemm5}
k_1 &=&\frac{V(t_{tr}) \Uc^{tr}}{ U(t_{tr})V(t_{tr})-AU(t_{tr}) + A \Uc^{tr}}\nonumber \\
k_2 &=& \frac{1}{\beta (A -V(t_{tr}))},
\end{eqnarray*}
being $\Uc^{tr} =\frac{c \delta}{\beta p(1-\eta_p)}$ the critical value for $U$ corresponding to the treatment parameter $p(1-\eta_p)$. 
Now, recalling that $\Rn(t_{tr}) = \frac{U(t_{tr})}{\Uc^{tr}}$, and denoting $U(t_{tr})$, $V(t_{tr})$ and $\Rn(t_{tr})$ as $U$, $V$ and $\Rn$, respectively, for the sake of simplicity, $k_1$, $k_2$ can be rewritten as:
\begin{eqnarray}\label{eq:lemm6}
k_1(\Rn) =\frac{ V}{ \frac{\delta}{\beta}( 1 -\Rn)^2 + \Rn V } > 0,~~~ \forall \Rn > 1
\end{eqnarray}
\vspace{-0.5cm}
\begin{eqnarray}\label{eq:lemm6_1}
k_2(\Rn) = \frac{1}{\delta (1-\Rn) - \beta V} < 0,~~~ \forall \Rn > 1
\end{eqnarray}
This way, $\hat t_{V,tr} = \hat t_{V,tr}(\Rn)=k_2(\Rn) \ln (k_1(\Rn)) + t_{tr}$ represents the time of the virus peak in terms of $\Rn$. Note that $\hat t_{V,tr}(\Rn)$ is defined only for $\Rn \geq 1$; indeed, $\hat t_{V,tr}(1) = t_{tr}$ since $k_1(1)=1$, while $\hat t_{V,tr}(\Rn) > t_{tr}$ for $\Rn > 1$.\\

The idea now is to consider the derivative of $\hat t_{V,tr}$ with respect to $\Rn$, to show that it is negative for treatment times $t_{tr}$ small enough. This derivative reads:
\begin{eqnarray*}\label{eq:lemm7}
	\frac{\partial \hat t_{V,tr}}{\partial \Rn} = \frac{\partial k_2 }{\partial \Rn} \ln(k_1) + k_2 \frac{1}{k_1} \frac{\partial k_1 }{\partial \Rn}
\end{eqnarray*}
where 
\begin{eqnarray}\label{eq:lemm8}
\frac{\partial k_1 }{\partial \Rn} =-~\frac{V (\frac{2\delta}{\beta} (\Rn-1) + V )  }{(\frac{\delta}{\beta} (1-\Rn)^2 + \Rn V)^2}
\end{eqnarray}
\vspace{-0.3cm}
\begin{eqnarray}\label{eq:lemm8_1}
\frac{\partial k_2 }{\partial \Rn} = \frac{\delta}{(\delta (1-\Rn) - \beta V)^2}.
\end{eqnarray}
Since $\frac{\partial k_2 }{\partial \Rn} >0$ for $\Rn > 1$, condition $\frac{\partial \hat t_{V,tr}}{\partial \Rn} <0$ can be written as
\begin{eqnarray}\label{eq:lemm9}
  \ln(k_1) <  - \frac{k_2}{k_1} \frac{\partial k_1 }{\partial \Rn} (\frac{\partial k_2 }{\partial \Rn})^{-1}.
\end{eqnarray}
Furthermore, given that $\ln(\cdot)$ is an increasing function on $(0,\infty]$, it follows that
\begin{eqnarray}\label{eq:lemm10}
k_1  < e^{- \frac{k_2}{k_1} \frac{\partial k_1 }{\partial \Rn} (\frac{\partial k_2 }{\partial \Rn})^{-1}}.
\end{eqnarray}

\begin{figure}[H]
    \centering
	\includegraphics[width=0.8\textwidth]{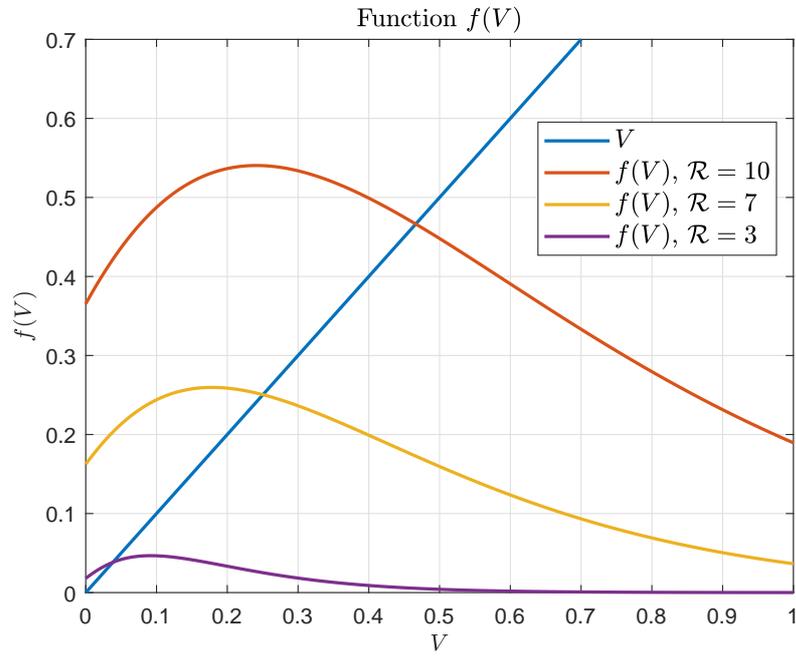}
	\caption{\small{Qualitative plot of function $f(V)$ (arbitrary parameters) for different values of $\mathcal{R}$. As it can be seen - independently of the parameter values - if $\Rn(t_{tr})>1$} it there exists an interval of values of $V$, and a corresponding period of time such that $f(V)>V$.}
	\label{fig:f_V}
\end{figure}

Then, by replacing \eqref{eq:lemm6}, \eqref{eq:lemm6_1}, \eqref{eq:lemm8} and \eqref{eq:lemm8_1} in inequality \eqref{eq:lemm10}, we have:
\begin{eqnarray*}\label{eq:lemm12}
\frac{\beta}{\delta} V <  [ ( 1 -\Rn)^2 + \Rn \frac{\beta}{\delta} V ] e^{ \frac{-2(1-\Rn)^2 + 2 \frac{\beta}{\delta} V (1-\Rn) - \frac{\beta}{\delta}^2V^2 }{(1-\Rn)^2 + \frac{\beta}{\delta} \Rn V}}
\end{eqnarray*}
Now, for $\Rn>1$ the exponent is negative, so function 
\begin{eqnarray*}\label{eq:lemm13}
	f(\frac{\beta}{\delta}V) :=  [ ( 1 -\Rn)^2 + \Rn \frac{\beta}{\delta} V ] e^{ \frac{-2(1-\Rn)^2 + 2 \frac{\beta}{\delta} V (1-\Rn) - \frac{\beta}{\delta}^2V^2 }{(1-\Rn)^2 + \frac{\beta}{\delta} \Rn V}}
\end{eqnarray*}
starts at the positive value $( 1 -\Rn)^2 e^{-2}$, for $V=0$, then increases to a maximum and finally decreases asymptotically to zero, for $V\rightarrow \infty$. So, there exists some interval of $V$, maybe small, such that $V<f(V)$. Figure \ref{fig:f_V} shows a plot of $f(V)$ for different values of $\Rn>1$.

\begin{figure}[H]
	\hspace{-1.5cm}
	\includegraphics[trim = 35mm 0mm 35mm 0mm, clip, width=1.2\textwidth]{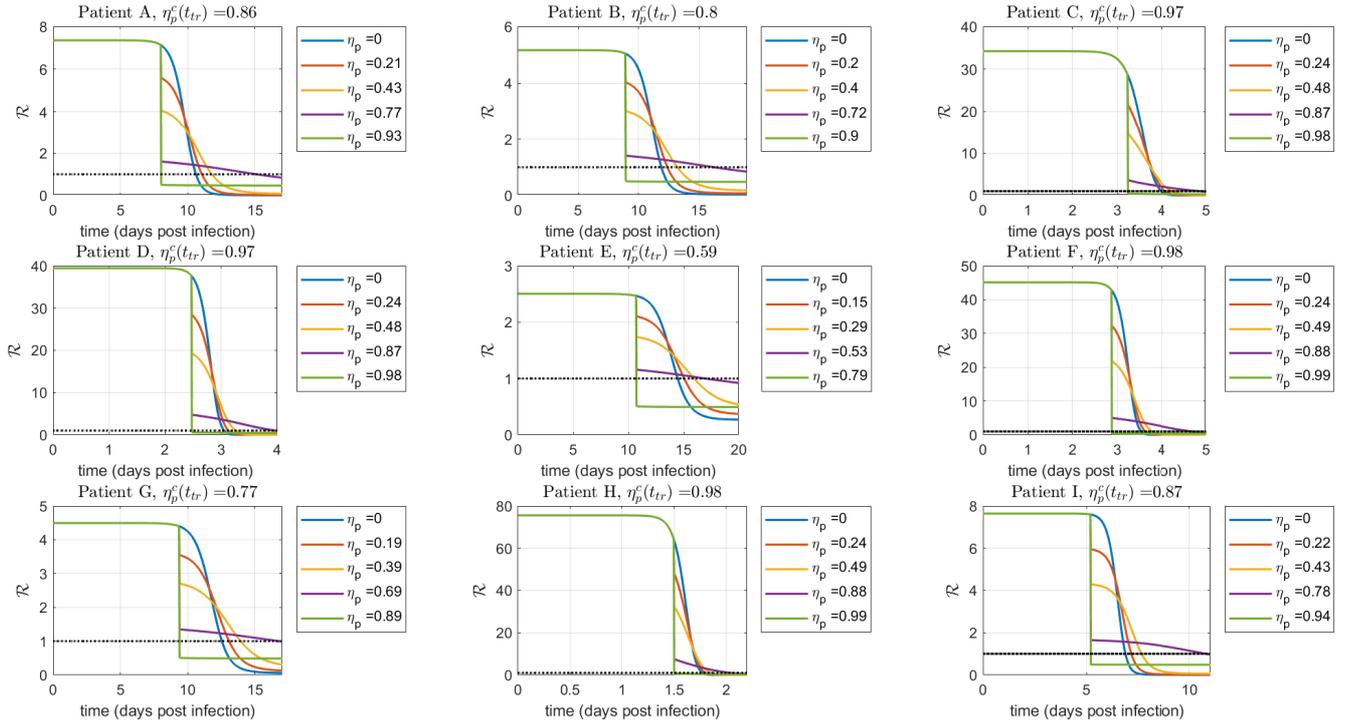}
	\caption{\small{Evolution of $\Rn(t)$ when an antiviral treatment is initiated at a time of approximately $0.75~ \hat t_V$, and different antiviral inhibition effects $\eta_p$ are considered, smaller and greater than the critical value $\eta_p^c$. The black dashed line represents $\Rn=1$. As it can be seen, for values of $\eta_p < \eta_p^c$, $\Rn(t)$ crosses $1$ at larger times for larger values of $\eta_p$ as it is stated in Lemma \ref{lem:Rt}. This implies that if $\eta_p < \eta_p^c$, larger values of $\eta_p$ delays the virus peak time, as it is stated in Theorem \ref{theo:mainresult}.iii}. Furthermore, the figure confirms that, for real patient date, $t^e$ is close to $\hat t_V$.}
	\label{fig:Rt}
\end{figure}

Finally, since we are considering the treatment time to belong to the increasing period of $V(t)$ (\textit{i.e.}, $t_{tr} \in (\check t_V,\hat t_V)$, with $V(\check t_V)\approx 0$, then small values of $V(t)$ correspond to small values of $t$. So, a time interval $(\check t_V,t^{e}]$ exists such that $\hat t_{V,tr}$ is a decreasing function of $\Rn(t_{tr})$, and the proof is complete. Figure \ref{fig:Rt} shows the time evolution of $\Rn(t)$ when the antiviral treatment is initiated at $0.75 \hat t_V$ and different antiviral inhibition effects are considered, for the real data patients simulated in Sections \ref{sec:pac} and \ref{sec:simres}. As it can be seen, larger values of $\eta_p$ (or smaller values of $\Rn(t_{tr})$) corresponds always with larger values of $\hat t_{tr,V}$. $\square$
\end{pf}

%----------------------------------------------------------------------------------------
%	BIBLIOGRAPHY
%----------------------------------------------------------------------------------------
\newpage
\bibliographystyle{elsarticle-num} 
\bibliography{biblioIR}

\end{document}